\documentclass[12pt]{article}

\input{xypic.tex} \xyoption{all}
\usepackage{amsthm,amssymb,latexsym}
\usepackage{graphics,color}

\newcommand{\real}{\mathbb{R}} \newcommand{\cpx}{\mathbb{C}}
\newcommand{\Z}{\mathbb{Z}} 
\newcommand{\p}{\mathbb{P}} \newcommand{\op}[1]{{\cal O}_{\mathbb{P}^{#1}}}
\newcommand{\bB}{\mathbb{B}} \newcommand{\vB}{\vec{B}}
\newcommand{\tW}{\tilde{W}}\newcommand{\cE}{{\cal E}}
\newcommand{\id}{\mathbf{1}}
\newcommand{\miq}{\mathfrak{M}^{\rm I}_q}
\newcommand{\mjq}{\mathfrak{M}^{\rm J}_q}
\newcommand{\mpq}{\mathfrak{M}_{p,q}}
\newcommand{\mijq}{\mathfrak{M}^{\rm IJ}_q}
 
\newcommand{\cQ}{{\cal Q}}
\newcommand{\inst}{{\cal M}^{\rm reg}(r,c)}     
\newcommand{\nak}{{\cal M}(r,c)}     
\newcommand{\cmc}{{\cal M}_{\cpx}(r,c)}        
\newcommand{\h}[1]{\otimes H^0(\mathbb{P}^1,{\cal O}_{\mathbb{P}^1}({#1}))}

\newcommand{\evp}{{\rm ev}_p} \newcommand{\del}{\partial}
\newcommand{\hf}{\frac{1}{2}}

\newtheorem{theorem}{Theorem}

\newtheorem{proposition}[theorem]{Proposition}

\newtheorem{corollary}[theorem]{Corollary}
\newtheorem{remark}[theorem]{Remark}

\newtheorem*{definition}{{\bf Definition}}

\begin{document}

\title{Complex ADHM equations, \\ sheaves on $\p^3$ \\ and quantum instantons}
\author{Igor B. Frenkel \\ Yale University \\ Department of Mathematics \\
10 Hillhouse Avenue \\ New Haven, CT 06520-8283 USA\\ \\
Marcos Jardim \\ IMECC - UNICAMP \\
Departamento de Matem\'atica \\ Caixa Postal 6065 \\
13083-970 Campinas-SP, Brazil}

\maketitle

\begin{abstract}
We use a complex version of the celebrated Atiyah-Hitchin-Drinfeld-Manin 
matrix equations to construct admissible torsion-free sheaves on $\p^3$
and complex quantum instantons over our quantum Minkowski space-time. We
identify the moduli spaces of various subclasses of sheaves on $\p^3$,
and prove their smoothness. We also define the Laplace equation in the
quantum Minkowski space-time, study its solutions and relate them to
the admissibility condition for sheaves on $\p^3$.
\end{abstract}

\newpage

\tableofcontents

\newpage
\baselineskip18pt


\section*{Introduction} \label{intro}

In \cite{FJ} we introduced a quantum Minkowski space-time based on the 
quantum group $SU(2)_q$ extended by a degree operator, and we formulated
a quantum version of the anti-self-dual Yang-Mills (ASDYM) equation with
unitary gauge group. A remarkable feature of the quantum equations is the
natural parameterization of their solutions of fixed rank $r$ and charge $c$
by the classical Atiyah-Hitchin-Drinfeld-Manin (ADHM) moduli spaces $\inst$.

The theory of ASDYM equations substantially involves the complex structures, 
including the fundamental twistor space of Penrose. It is therefore appropriate
to study also the ASDYM equations with a complex linear gauge group on the
complexified Minkowski space-time (see e.g. \cite{WW}).
The corresponding counterparts of the ADHM moduli spaces $\inst$ were 
implicitly described by Donaldson in \cite{D1}, who in fact used them to
construct certain holomorphic vector bundles on $\p^3$ related to instantons
by the Penrose twistor diagram. These holomorphic vector bundles
admit a cohomological characterization \cite{Ma}, and are known in the 
literature as complex instanton bundles \cite{OSS}.

In this paper, we expand the classical theory of the complex ASDYM equations
in two related directions. On one hand, we will consider the moduli spaces of
complex {\em torsion free instanton sheaves} on $\p^3$ and show that they
correspond to extending the moduli spaces of complex instantons
${\cal M}_{\cpx}^{\rm reg}(r,c)$ to the larger moduli spaces of stable complex
ADHM data, denoted by $\cmc$. On the other hand, we will use the stable complex
ADHM data to construct {\em complex quantum instantons} over the complexified
quantum Minkowski space-time introduced in \cite{FJ}. In this paper we also
generalize the classical relation \cite{Ma} between the admissibility
condition that determines the instanton bundles and the solutions of the Laplace
equations by developing the quantum differential calculus.

The organization of the paper is completely reflected in the title, and
consists of three sections complemented by three appendices. In Section 1
we define the complex ADHM equations and the associated moduli spaces.
We prove various results about them, including smoothness, dimension and
non-emptiness of $\cmc$ and its open subset ${\cal M}_{\cpx}^{\rm reg}(r,c)$
for $r\geq2$. In Section 2, we define and study admissible torsion-free
sheaves on $\p^3$. The main result of this section is the characterization
of the sheaves corresponding to $\cmc$. Using such characterization we show
that the spaces ${\cal M}_{\cpx}(1,c)$ are empty. We also identify the open
subsets of $\cmc$ that characterize reflexive and locally-free admissible
sheaves. In particular, we verify that the moduli spaces of framed complex
instanton bundles of rank $r$ and charge $c$ over $\p^3$ is given by the
regular complex ADHM data ${\cal M}_{\cpx}^{\rm reg}(r,c)$ and has the
structure of a smooth complex manifold of dimension $4rc$. In Section 3
we recall the definition of the quantum Minkowski space-time and the
quantum ASDYM equations, and we construct quantum instantons from the
ADHM data parameterized by the moduli spaces $\cmc$. At the last subsection
we define the quantum Laplace equation and produce their solutions
from the cohomology of an instanton sheaf on $\p^3$. This extends the phenomenon
observed in our first paper \cite{FJ}, in the case of anti-self-dual Yang-Mills
equation: the solutions of the {\em quantum equation} are parameterized by the
{\em classical data}. We thus obtain a surprisingly
{\em explicit and direct link between the commutative geometry encoded into sheaves
on $\p^3$ and the noncommutative geometry of the quantum Minkowski space-time}.
We also discuss the consistency of the Laplace equations on the two affine
parts $\miq$ and $\mjq$ of our quantum Minkowski space-time, thus relating
the admissibility condition with the non-existence of consistent solutions
to the quantum Laplace equation. Finally, in the three appendices corresponding
to each of the three sections we collect various
background facts and calculations used throughout the paper.

Our constructions are best reflected by the following triangle, whose 
vertices correspond to the objects studied in the three sections:
\vskip18pt 
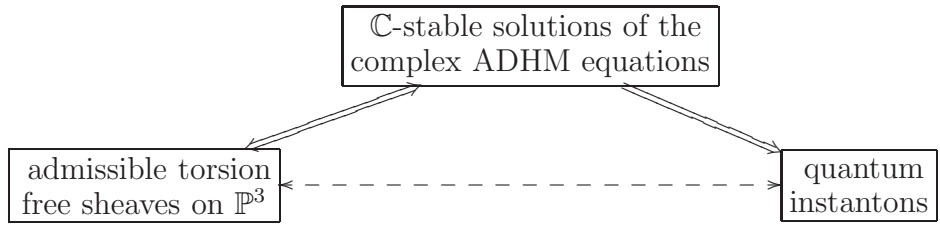
\centerline{ \xymatrix{
& *+[F]\txt{ $\cpx$-stable solutions of the \\ complex ADHM equations }
\ar@{<=>}[dl] \ar@{=>}[dr] & \\
*+[F]\txt{ admissible torsion \\ free sheaves on $\p^3$ }
\ar@{<-->}[rr] & &
*+[F]\txt{ quantum \\ instantons }
}} \vskip18pt
The solid arrows were established in this paper; the dashed arrow corresponds
to the quantum Penrose transform introduced in \cite{FJ} and which we plan
to study further in a sequel to this paper. The completion of this circle of
ideas will yield a complete characterization of quantum instantons, while
opening, at the same time, a new perspective in the theory of sheaves on
$\p^3$. It should also enhance the direct connection of the commutative
geometry encoded into sheaves on $\p^3$ and the noncommutative geometry
of the quantum Minkowski space-time.

\bigskip

\paragraph{Acknowledgment.}
I.F. is supported by the NSF grant DMS-0070551.
M.J. thanks Eyal Markman for useful discussions.


\section{Complex ADHM equations} \label{rev}

We start with some fundamental definitions and further motivation,
before formulating our first main result in Section \ref{adhm}.

\subsection{ADHM data and instantons}
Let $V$ and $W$ be complex vector spaces, with dimensions
$c$ and $r$, respectively, and consider maps $B_1,B_2\in{\rm End}(V)$,
$i\in{\rm Hom}(W,V)$ and $j\in{\rm Hom}(V,W)$. This so-called
{\em ADHM datum} $(B_1,B_2,i,j)$ is said to be:
\begin{enumerate}
\item {\em stable}, if there is no proper subspace $S\subset V$ such that
$B_k(S)\subset S$ ($k=1,2$) and $i(W)\subset S$;
\item {\em costable}, if there is no proper subspace $S\subset V$ such that
$B_k(S)\subset S$ ($k=1,2$) and $S\subset \ker j$;
\item{\em regular}, if it is both stable and costable.
\end{enumerate}
Notice that $(B_1,B_2,i)$ is stable if and only if the triple 
$(B_1^*,B_2^*,i^*)$ consisting of the dual maps is costable.

Now, providing $V$ and $W$ with Hermitian structures, we consider the
so-called {\em ADHM equations} ($\dagger$ denotes Hermitian conjugation):
\begin{eqnarray}
\label{adhm1} [ B_1 , B_2 ] + ij & = & 0  \\
\label{adhm2} [ B_1 , B_1^\dagger ] + [ B_2 , B_2^\dagger ] + ii^\dagger -
j^\dagger j & = & \xi \id_V
\end{eqnarray}
With $GL(V)$ acting on the set of all ADHM data in the following way:
$$ g(B_{kl},i_k,j_k)=(gB_{kl}g^{-1},gi_k,j_lg^{-1}),
~~~ g\in GL(V) $$
we define the following varieties:
$$ {\cal M}^0(r,c) = \left. \left\{ \rm{all ~ solutions ~ of} ~
(\ref{adhm1}) ~ \rm{and} ~ (\ref{adhm2}) ~ {\rm with} ~ \xi=0 \right\} 
\right/ U(V) $$
\begin{eqnarray*}
\nak & = & \left. \left\{ \rm{stable ~ solutions ~ of} ~
(\ref{adhm1}) \right\} \right/ GL(V) ~ \simeq ~ \\
& \simeq & \left. \left\{ \rm{all ~ solutions ~ of} ~ (\ref{adhm1}) ~ 
\rm{and} ~ (\ref{adhm2})
~ {\rm with} ~ \xi>0 \right\} \right/ U(V)
\end{eqnarray*}
\begin{eqnarray*}
\inst & = & \left. \left\{ \rm{regular ~ solutions ~ of} ~
(\ref{adhm1}) \right\} \right/ GL(V) ~ \simeq ~ \\
& \simeq & \left. \left\{ \rm{regular ~ solutions ~ of} ~
(\ref{adhm1}) ~ \rm{and} ~ (\ref{adhm2}) ~ {\rm with} ~ \xi=0 \right\} 
\right/ U(V)
\end{eqnarray*}
It can be shown that $\inst$ and $\nak$ are smooth hyperk\"ahler manifolds
of dimension $2rc$, with the exception that $\inst$ is empty for $r=1$; the
proofs are given in Appendix \ref{ap-nak} for they are used in the next section.
Moreover, $\inst$ is the smooth locus of the singular variety ${\cal M}^0(r,c)$
for $r\geq2$, while ${\cal M}(r,c)$ is the minimal resolution of ${\cal M}^0(r,c)$
\cite{N1,N2}.


\subsection{Complex ADHM data} \label{adhm}

As above, let $V$ and $W$ be complex vector spaces, with dimensions
$c$ and $r$, respectively; Hermitian structures are no longer required.
Set $\tW=V\oplus V\oplus W$, and define also:
$$ \mathbf{B} = {\rm Hom}(V,V)\oplus{\rm Hom}(V,V)\oplus
{\rm Hom}(W,V)\oplus{\rm Hom}(V,W) $$
Let $\bB=\mathbf{B}\oplus\mathbf{B}$, with $\vB=(B_{kl},i_k,j_k)\in\bB$
($k,l=1,2$) being called a {\it complex ADHM datum}. As usual, the group 
$GL(V)$ acts naturally on $\mathbf{B}$ and on $\bB$, in the following way:
\begin{equation} \label{action}
g(B_{kl},i_k,j_k)=(gB_{kl}g^{-1},gi_k,j_kg^{-1}),
\ \ \ g\in GL(V)
\end{equation}
Equivalently, we can think of an element in $\bB$ as a holomorphic
section of the bundle $\mathbf{B}\otimes\op1(1)$ by defining:
\begin{equation} \label{st1}
\tilde{B}_1 = zB_{11} + wB_{21} \ \ \ {\rm and} \ \ \
\tilde{B}_2 = zB_{12} + wB_{22}
\end{equation}
\begin{equation} \label{st2}
\tilde{\imath} = zi_1 + wi_2 \ \ \ {\rm and} \ \ \
\tilde{\jmath} = zj_1 + wj_2
\end{equation}
In other words, $\bB=\mathbf{B}\h1 $, with $z,w$ denoting a basis of
$H^0(\p^1,\op1(1))$ or, equivalently, a choice of homogeneous
coordinates in $\p^1$. In particular, one can also view the maps
(\ref{st1}) and (\ref{st2}) in the following way:
$$ \tilde{B}_1,\tilde{B}_2 \in {\rm Hom}(V,V)\h1 $$
$$ \tilde{\imath} \in {\rm Hom}(W,V)\h1 ~~ {\rm and} ~~
\tilde{\jmath} \in {\rm Hom}(V,W)\h1 $$

The evaluation map $\evp:H^0(\p^1,\op1(1))\to\cpx$ can be tensored with
the identity to yield maps $\evp:\bB\to\mathbf{B}$ and
$\evp:{\rm Hom}(V,V)\h1\to{\rm Hom}(V,V)$. For simplicity, we use the
notation $\vB_p=\evp(\vB)$.

\begin{definition} \label{st-defn}
A complex ADHM datum $\vB=(B_{kl},i_k,j_k)$ is said to be:
\begin{enumerate}
\item $\cpx$-{\em semistable} if there is $p\in\p^1$ such that $\vB_p$ is 
stable;
\item $\cpx$-{\em stable} if $\vB_p$ is stable for all $p\in\p^1$;
\item $\cpx$-{\em costable} if $\vB_p$ is costable for all $p\in\p^1$;
\item $\cpx$-{\em semiregular} if it is $\cpx$-stable and there is 
$p\in\p^1$
such that $\vB_p$ is regular;
\item $\cpx$-{\em regular} if $\vB_p$ is regular for all $p\in\p^1$.
\end{enumerate} \end{definition}

The motivation behind these definitions will be clearer in the next Section:
$\cpx$-stable, $\cpx$-semiregular, and $\cpx$-regular will correspond to
torsion-free, reflexive and locally-free sheaves on $\p^3$, respectively.
In particular, notice that $\vB=(B_{kl},i_k,j_k)$ is $\cpx$-regular if and 
only if it is both $\cpx$-stable and $\cpx$-costable. Moreover, if 
$\vB=(B_{kl},i_k,j_k)$ is $\cpx$-stable then $(B_{11},B_{12},i_1,j_1)$ and
$(B_{21},B_{22},i_2,j_2)$ are both stable.

\begin{proposition} \label{free}
If $\vB$ is $\cpx$-semistable, then its $GL(V)$ stabilizer is trivial.
\end{proposition}
\begin{proof}
If $\vB$ is fixed by some nontrivial $g\in GL(V)$, then $\vB_p$ is also
fixed by $g$ for all $p\in\p^1$. Thus, by Proposition \ref{stability} in
Appendix \ref{ap-nak}, there is no $p\in\p^1$ such that $\vB_p$ is stable.
\end{proof}

The first main goal of this paper is to study the {\em complex ADHM equations}:
\begin{eqnarray}
\label{c1} [ B_{11} , B_{12} ] + i_1j_1 & = & 0 \\
\label{c2} [ B_{21} , B_{22} ] + i_2j_2 & = & 0 \\
\label{c3} [ B_{11} , B_{22} ] + [ B_{21} , B_{12} ] +
i_1j_2 + i_2j_1 & = & 0
\end{eqnarray}
which were first posed by Donaldson in \cite{D1}; it is important to note 
that the equations (\ref{c1}-\ref{c3}) are equivalent to:
\begin{equation} \label{c4}
[ \tilde{B}_1 , \tilde{B}_2 ] + \tilde{\imath}\tilde{\jmath} = 0,\ \ \forall 
[z:w] \in \p^1
\end{equation}

It is easy to see that solutions of (\ref{c1}-\ref{c3}) are preserved by the
$GL(V)$ action (\ref{action}). Therefore, we define the moduli space of 
solutions of the complex ADHM equations as the quotient:
$$ \cmc := \left. \left\{
\begin{array}{c} \cpx-{\rm stable} \\
{\rm solutions\ of\ (\ref{c1}-\ref{c3})} \end{array}
\right\} \right/ GL(V) $$
Our first main result, to be proved in this Section, states that $\cmc$ is a
smooth, complex manifold of complex dimension $4rc$, non-empty for $r\geq2$.
The strategy of the proof is the same as for Theorem \ref{nakvar} in Appendix
\ref{ap-nak} ; we consider the map
\begin{eqnarray*}
& \tilde{\mu} : \bB^{\rm st}\to{\rm Hom}(V,V)\h2 & \\
& \tilde{\mu} (B_{kl},i_k,j_k) = [ \tilde{B}_1 , \tilde{B}_2 ] + 
\tilde{\imath}\tilde{\jmath} &
\end{eqnarray*}
where $\bB^{\rm st}$ denote the open subset of $\cpx$-stable complex ADHM 
data. Clearly $\cmc=\tilde{\mu}^{-1}(0)/GL(V)$. We have already established that
$GL(V)$ acts freely on $\bB^{\rm st}$; we must then argue that $\tilde{\mu}$
has a surjective derivative and that the $GL(V)$ action has a closed graph. 

\begin{proposition} \label{surjectivity}
$\vB$ is $\cpx$-stable if and only if the derivative map
$$ D_{\vB}\tilde{\mu}:\bB\to{\rm Hom}(V,V)\h2 $$
is surjective.
\end{proposition}

\begin{proof}
Considering the maps ($k=1,2$):
\begin{eqnarray*}
& \partial_k : \mathbf{B}\to{\rm Hom}(V,V) & \\
& \partial_k = [ \cdot , B_{k2} ] + [ B_{k1} , \cdot ] + i_k\cdot + \cdot 
j_k &
\end{eqnarray*}
which can also be regarded as $c^2\times(2c^2+2cr)$ matrix. We can then 
express the $3c^2\times(4c^2+4cr)$
matrix of the derivative map $D_{\vB}\tilde{\mu}$ in the following form:
\begin{equation} \label{matrix} D_{\vB}\tilde{\mu} = \left(
\begin{array}{cc} \partial_1 & 0 \\ \partial_2 & \partial_1 \\ 0 & 
\partial_2 \end{array}
\right) \end{equation}
Our goal is to show that the above matrix has maximal rank $3c^2$ if and 
only if $\vB$ is $\cpx$-stable.

So let $\{l_k\}_{k=1}^{c^2}$ denote the rows of the matrix $\partial_1$, and 
let $\{l'_k\}_{k=1}^{c^2}$ denote the rows of the matrix $\partial_2$; each
$l_k,l'_k$ is regarded as a vector in $\mathbf{B}$. As remarked above,
$\cpx$-stability is equivalent to the vectors $\{zl_k+wl'_k\}_{k=1}^{c^2}$
being linearly independent (as vectors in $\mathbf{B}$) for all $[z:w]\in\p^1$. 
The rows of the matrix $(\partial_2 ~~ \partial_1)$ are then given by
$\{(l'_k,l_k)\}_{k=1}^{c^2}$ regarded as vectors in $\mathbf{B}\oplus\mathbf{B}$.
Clearly,  the matrix (\ref{matrix}) above is surjective if and only if
$\{(l_k,0)~,~(l'_k,l_k)~,~(0,l'_k)\}_{k=1}^{c^2}$ form a linearly 
independent set of vectors in $\bB$; this is in turns equivalent to the
statement that if the coefficients $\gamma^k$ are such that
$\sum_k\gamma^kl_k\in{\rm span}\{l'_k\}$ and
$\sum_k\gamma^kl'_k\in{\rm span}\{l_k\}$, then $\gamma^k=0$.

Let $L={\rm span}\{l_k\}_{k=1}^{c^2}$, $L'={\rm span}\{l'_k\}_{k=1}^{c^2}$; 
the theorem can then be reduced to the following statement: the vectors 
$\{zl_k+wl'_k\}_{k=1}^{c^2}$ are linearly independent
for all $[z:w]\in\p^1$ if and only if $\sum_k\gamma^kl_k\in L'$ and 
$\sum_k\gamma^kl'_k\in L$ implies $\gamma^k=0$.

First, assume that if $\sum_k\gamma^kl_k\in L'$ and
$\sum_k\gamma^kl'_k\in L$, then $\gamma^k=0$. If
$\sum_k \gamma^k(zl_k+wl'_k)=0$, then
$L\ni\sum_k\gamma^kl_k=-\frac{w}{z}\sum_k\gamma^kl'_k\in L'$;
hence $\gamma^k=0$, and $\{zl_k+wl'_k\}_{k=1}^{c^2}$ is linearly
independent for all $[z:w]\in\p^1$.

For the converse direction, denote $I=L\cap L'$; we can assume that
$I={\rm span}\{l_k\}_{k=1}^d={\rm span}\{l'_k\}_{k=1}^d$. Since 
$\sum_k\gamma^kl_k,~\sum_k\gamma^kl'_k\in I$,
we have $\gamma^k=0$ for $k=d+1,\dots,c^2$. Moreover, for each $k=1,\dots,d$, 
we have $l'_k=\sum_{n=1}^d g^n_kl_n$; let $G$ be the invertible $d\times d$
matrix with entries $g^n_k$. Therefore,
$$ \sum_{k=1}^d \gamma^k(zl_k+wl'_k) = \sum_{k,n=1}^d \gamma^k
\left( z\delta_k^n+wg^n_k  \right) l_n =  \sum_{n=1}^d c^n(z,w)l_n ~ ,$$
where $c^n(z,w)$ are the entries of the vector $\gamma(z\id+wG)$; note that, 
by construction, the matrix $(z\id+wG)$ is also invertible for generic $[z:w]\in\p^1$.
Now, the vectors $\{l_n\}$ are linearly independent, so if
$\sum_{k=1}^d \gamma^k(zl_k+wl'_k) = 0$, then $c_n(z,w)=0$ for each
$n=1,\dots,d$ and for all $[z:w]\in\p^1$. This implies that $\gamma^k=0$,
as desired.
\end{proof}

As a by-product of our proof, we obtain the following interesting result:

\begin{proposition} \label{ss}
If the complex ADHM datum $\vB$ is:
\begin{itemize}
\item  $\cpx$-semistable, then there are at most finitely many
points $p\in\p^1$ such that $\vB_p$ is not stable.
\item $\cpx$-semiregular, then there are at most finitely many
points $p\in\p^1$ such that $\vB_p$ is not regular.
\end{itemize} \end{proposition}
\begin{proof}
Assume that $\vB_p$ is stable for some $p=[p_1:p_2]\in\p^1$. Then 
$p_1\partial_1+p_2\partial_2$ is surjective and its row vectors
$\{p_1l_k+p_2l'_k\}_{k=1}^{c^2}$ are linearly independent as
vectors in $\mathbf{B}$. We complete $\{p_1l_k+p_2l'_k\}_{k=1}^{c^2}$ to a 
basis of $\mathbf{B}$, and denote by $H(p)$ the $(2c^2+2cr)\times(2c^2+2cr)$
matrix formed by such basis.

As we vary $p\in\p^1$, we see that $\det H(p)$ is a section of $\op1(c^2)$. 
Hence $\{zl_k+wl'_k\}_{k=1}^{c^2}$ must be linearly independent except for 
finitely many $[z:w]\in\p^1$, which means that $\vB_{[z:w]}$ is stable away
from finitely many (up to $c^2$) points in $\p^1$.

The second statement follows by duality.
\end{proof}

\begin{proposition}
The action (\ref{action}) has a closed graph.
\end{proposition}
\begin{proof}
Let $\{X_k\}$ be a sequence in $\bB^{\rm st}$, while $\{g_k\}$ denotes a
sequence in $GL(V)$; assuming that
$$ \lim_{k\to\infty} X_k=X_\infty ~~~ {\rm and} ~~~
\lim_{k\to\infty} g_k\cdot X_k=Y_\infty ~,$$
we show that the sequence $\{g_k\}$ converges to some $g_\infty\in GL(V)$,
so that $Y_\infty =g_\infty\cdot X_\infty$.

Using evaluation at $p\in\p^1$, it follows that:
$$ \lim_{k\to\infty} (X_k)_p = (X_\infty)_p ~~~ {\rm and} ~~~
\lim_{k\to\infty} g_k \cdot (X_k)_p = (Y_\infty)_p ~ . $$
Hence, by argument in the proof of Proposition \ref{proper} in Appendix
\ref{ap-nak}, we conclude that
$$ \lim_{k\to\infty} g_k = g_\infty =
R\left( (Y_\infty)_p \right) (T_\infty)_p \left[ R\left( (X_\infty)_p \right) (T_\infty)_p \right]^{-1} \in GL(V) ~,$$
where the map $R (X) : W^{\oplus c^2} \longrightarrow V$ given, for
$X=(B_1,B_2,i,j)\in\mathbf{B}^{\rm st}$, by:
$$ R(X) = i \oplus \cdots \oplus B_1^mB_2^n i \oplus \cdots \oplus B_1^{c-1}B_2^{c-1} i
~~,~~ 1\leq m,n \leq c-1 ~ .$$
Note that $gR(X)=R(g\cdot X)$ for any $g\in GL(V)$, and that 
$R(X)$ is surjective if and only if $X$ is stable.
\end{proof}


\subsection{Existence of solutions}

So far, we can conclude from Propositions \ref{free} and \ref{surjectivity} 
that $\cmc$ is a smooth complex manifold of dimension $4rc$ provided it is 
non-empty. The case $r=\dim W=1$ is rather special due to the following result:

\begin{proposition} \label{no}
There are no $\cpx$-stable solutions of (\ref{c1}-\ref{c3}) for $r=1$.
\end{proposition}

The proof will be delayed until the end of Section \ref{adhm2p3} (see 
Proposition \ref{no.r=1}). However, let us consider here the simplest
possible case: $r=c=1$. In this case, a $\cpx$-stable solution of
(\ref{c1}-\ref{c3}) reduces to six complex numbers $(b_{kl},i_k)$,
since $j_1=j_2=0$ by Proposition \ref{r=1case}. Now $\tilde{\imath}$
is simply a section of $\op1(1)$, so it must vanish at one point
$p\in\p^1$, which implies that $\evp(\vB)$ is not stable.

This example also illustrate the fact that {\em there exist 
$\cpx$-semistable solutions of (\ref{c1}-\ref{c3}) which are
not $\cpx$-stable}.

Fortunately, the existence of regular solutions for the real ADHM equations 
(\ref{adhm1}) and (\ref{adhm2}) for $r\geq 2$ can be used to guarantee the
existence of $\cpx$-stable solutions of the complex equations.

Indeed, note that if $V$ are $W$ are provided with a Hermitian inner 
product, then the space of complex ADHM data $\bB$ acquires a natural
involution $\dagger:\bB\to\bB$ given by:
$$ \dagger(B_{11},B_{12},B_{21},B_{22},i_1,i_2,j_1,j_2)=
(B^\dagger_{22},-B^\dagger_{21},-B^\dagger_{12},B^\dagger_{11},
j_2^\dagger,-j_1^\dagger,-i_2^\dagger,i_1^\dagger) $$
The point $\vB\in\bB$ is said to be {\em real} if it is fixed by $\dagger$.

Note that if $\vB$ is real, then complex ADHM equations (\ref{c1}-\ref{c3})
above reduce to the usual ADHM equations with $\xi=0$, by setting
$B_1=B_{11}=B_{22}^\dagger$, $B_2=B_{12}=-B_{21}^\dagger$,
$i=i_1=j_2^\dagger$ and $j=j_1=-i_2^\dagger$.

\begin{proposition} \label{c-st=reg}
If $(B_1,B_2,i,j)$ is a stable (hence regular) solution of (\ref{adhm1}) and 
(\ref{adhm2}) with $\xi=0$,
then $\vB=(B_1,B_2,-B_2^\dagger,B_1^\dagger,i,-j^\dagger,j,i^\dagger)$ is a 
$\cpx$-regular solution of
(\ref{c1}-\ref{c3}).
\end{proposition}

In particular, it follows that $\cmc$ is non-empty for all $r\geq2$ and for 
all $c\geq1$. Moreover, there is a holomorphic surjective map:
$$ \varepsilon : {\cal M}^{\rm reg}_{\cpx}(r,c) \to \inst ~~ , ~~ \forall 
r\geq2, k\geq1 $$
$$ \varepsilon (B_{kl},i_k,j_k) := (B_{11},B_{12},i_1,j_1) $$
Clearly, $\varepsilon$ is surjective, and the fibers
$\varepsilon^{-1}(B_1,B_2,i,j)$ are closed subsets of $\inst$ of
dimension $2rc$.

\begin{proof}
It is easy to see that $(B_1,B_2,i,j)$ satisfies (\ref{adhm1}) and 
(\ref{adhm2}) with $\xi=0$, if and only if $\vB$ as above satisfies
(\ref{c1}-\ref{c3}). Now note that in this case:
$$ \tilde{B}_1= zB_1-wB_2^\dagger ~~,~~ \tilde{B}_2= zB_2+wB_1^\dagger ~~,~~ 
\tilde{\imath}=zi-wj^\dagger ~.$$
If $\vB$ is not $\cpx$-stable, there is $[z:w]\in\p^1$ and a proper subspace 
$S\subset V$ such that
$[\tilde{B}_1^\dagger,\tilde{B}_2^\dagger]|_{S} = 0$ and 
$S\subset\ker~\tilde{\imath}^\dagger$. Thus
$i^\dagger|_{S} = k\cdot j|_{S}$ for some $k\in\cpx$, hence $ii^\dagger|_{S} 
= k\cdot ij|_{S}=[B_1,B_2]|_{S}$.
Hence ${\rm Tr}(ii^\dagger|_{S})=0$, so that $[B_1^\dagger,B_2^\dagger]|_{S} 
= 0$ and $S\subset\ker~i^\dagger$
and $(B_1,B_2,i,j)$ is not stable.

Thus we conclude that $\vB$ as above is $\cpx$-stable. However, it is not 
difficult to see that every real, $\cpx$-stable complex ADHM datum is
$\cpx$-regular. Indeed, if $\vB$ is real, then:
$$ \begin{array}{lr}
\tilde{B}_1 = zB_{11} - wB_{12}^\dagger & \tilde{B}_2 = zB_{12} + 
wB_{11}^\dagger \\
\tilde{\imath} = zi_1 - wj_1^\dagger & \tilde{\jmath} = z j_1 + wi_1^\dagger
\end{array} $$
Thus if $(\tilde{B}_1,\tilde{B}_2,\tilde{\jmath})$ is not costable at 
$[z:w]$, then $(\tilde{B}_1,\tilde{B}_2,\tilde{\imath})$ is not stable at
$[-\overline{w}:\overline{z}]$.
\end{proof}

\begin{remark} \label{st.not.reg} \rm
It is interesting to note that, differently from the real ADHM equations, 
{\em there are $\cpx$-stable solutions of (\ref{c1}-\ref{c3}) which are not 
$\cpx$-semiregular} (compare with Proposition \ref{r/ir}). Indeed, for $r=2$
and $c=1$, we can take:
$$ B_{kl}=0 ~~,~~ i_1=\left(\begin{array}{c} 1 \\ 0 \end{array}\right) ~~,~~
i_2=\left(\begin{array}{c} 0 \\ 1 \end{array}\right) ~~,~~ j_1=j_2=0 ~ . $$
Furthermore, {\em there are $\cpx$-semiregular solutions of (\ref{c1}-\ref{c3})
which are not $\cpx$-regular}; for $r=3$ and $c=1$, we can take: 
$$ B_{kl}=0 ~~,~~ i_1=\left(\begin{array}{c} 1 \\ 0 \\ 0 \end{array}\right)
~~,~~ i_2=\left(\begin{array}{c} 0 \\ 1 \\ 0 \end{array}\right) ~~,~~
j_1=(0 ~ 0 ~ 1) ~~,~~ j_2=0 ~ . $$
However, as we shall see in Section 3, every $\cpx$-semiregular solution of
(\ref{c1}-\ref{c3}) for $r=2$ is in fact $\cpx$-regular.
\end{remark}

Our first man result follows from Proposition \ref{c-st=reg} together with
Propositions \ref{free} and \ref{surjectivity}:

\begin{theorem} \label{main1}
$\cmc$ is a smooth complex manifold of dimension $4rc$, non-empty for 
$r\geq2$, $c\geq1$.
\end{theorem}

Given the close analogy between $\cmc$ and the varieties $\inst$ and $\nak$,
we expect $\cmc$ to be a connected, hyperk\"ahler quasi-projective algebraic
variety for all $r\geq2$ and $c\geq1$. Indeed, one has for each $p\in\p^1$ a
holomorphic map $\varepsilon_p:\cmc\to\nak$ induced by the natural evaluation
map $\evp:\bB\to\mathbf{B}$ introduced above. One can then try to establish
various properties of $\cmc$ by studying various properties (e.g. surjectivity)
of the evaluation map.

\subsection{Solutions for $c=1$}

For $c=\dim~V=1$, the varieties $\cmc$ can be described quite concretely.
In this case, $B_{kl}$ are just complex numbers, while $i_k$ and $j_k$
can be regarded as vectors in $W$; the complex ADHM equations reduce to:
\begin{equation} \label{adhmc.c=1}
i_1j_1 = i_2j_2 = i_1j_2 + i_2j_1 = 0 ~ .
\end{equation}
It is then easy to see that ${\cal M}_{\cpx}(r,1)=\cpx^4\times{\cal B}(r)$, 
where ${\cal B}(r)$ is a smooth quasi-projective variety of dimension
$4(r-1)$, described as follows.

Setting:
$$ i_1 = (x_1 ~ \cdots ~ x_r) ~~~  i_2 = (y_1 ~ \cdots ~ y_r) $$
$$ j_1 =\left( \begin{array}{c} z_1 \\ \vdots \\ z_r \end{array} \right) ~~~
j_2 =\left( \begin{array}{c} w_1 \\ \vdots \\ w_r \end{array} \right) $$
the equations (\ref{adhmc.c=1}) become:
\begin{equation} \label{bquadrics}
\sum_{k=1}^r x_kz_k = \sum_{k=1}^r y_kw_k =
\sum_{k=1}^r x_kw_k + y_kz_k  =  0 ~ .
\end{equation}
while $\cpx$-stability is equivalent to the vectors $(x_1,\dots,x_r)$ and 
$(y_1,\dots,y_r)$ being linearly independent.

Then ${\cal B}(r)$ is the complete intersection of the three quadrics 
(\ref{bquadrics}) within the open subset of $\p^{4r-1}=\p(W^{\oplus4})$
consisting of $\cpx$-stable points.
In particular, we conclude that ${\cal M}_{\cpx}(r,1)$ is quasi-projective 
for all $r\geq 2$, in partial support of our general conjecture.



\section{Sheaves on $\mathbf{\p^3}$} \label{p3}

In this section we will characterize $\cmc$ as a moduli space of certain 
sheaves on $\p^3$. First, we recall the following definition, due to Manin
\cite{Ma}.

\begin{definition}
A coherent sheaf $\cE$ on $\p^3$ is said to be {\em admissible} if
$H^p(\p^3,\cE(k)) = 0$ for $p\leq1$, $p+k\leq-1$ and for $p\geq2$, 
$p+k\geq0$.
\end{definition}

This somewhat mysterious cohomological condition is made natural once we
recall that, under Penrose transform, locally-free admissible sheaves on
$\p^3$ correspond to (framed) $GL(r,\cpx)$ instantons (see Section
\ref{4.3} below). With a few extra assumptions on $\cE$, the admissibility
condition becomes a lot simpler; let $\ell_\infty=\{z=w=0\}$.


\begin{proposition} \label{fund.lemma}
Let $\cE$ be a torsion-free sheaf on $\p^3$ such that
$\cE|_{\ell_\infty}\simeq{\cal O}_{\ell_\infty}^{\oplus r}$.
$\cE$ is admissible if and only if $H^1(\p^3,\cE(-2)) = H^2(\p^3,\cE(-2)) = 0$.
\end{proposition}
\begin{proof}
Let $\wp$ be a plane containing $\ell_\infty$, e.g. $\wp=\{z=0\}$. Then
$\cE|_{\wp}$ is a torsion-free sheaf on $\wp$ which is trivial at $\ell_\infty$.
From \cite{N2} we know that:
\begin{equation} \label{v1}
H^0(\wp,\cE|_{\wp}(k)) = 0 ~ \forall k\leq-1 ~~ {\rm and} ~~
H^2(\wp,\cE|_{\wp}(k)) = 0 ~ \forall k\geq-2 ~~.
\end{equation}
Now consider the sheaf sequence:
\begin{equation} \label{v2}
0 \to \cE(k-1) \stackrel{\cdot z}{\longrightarrow} \cE(k)
\longrightarrow \cE|_{\wp}(k) \to 0
\end{equation}
Using (\ref{v1}), we conclude that:
$$ H^3(\p^3,\cE(k)) =  H^3(\p^3,\cE(k-1)) ~ \forall k\geq-2 $$
But, by Serre's vanishing theorem, $H^3(\p^3,\cE(N))=0$ for sufficiently 
large $N$, thus
$H^3(\p^3,\cE(k)) = 0$ for all $k\geq-3$.

Similarly, we have:
$$ H^0(\p^3,\cE(k-1)) =  H^0(\p^3,\cE(k)) ~ \forall k\leq-1 $$
Since $\cE\hookrightarrow \cE^{**}$, we have via Serre duality:
$$ H^0(\p^3,\cE(k))\hookrightarrow H^0(\p^3,\cE^{**}(k)) = 
H^3(\p^3,\cE^{***}(-k-4))^*~. $$
Thus, again by Serre's vanishing theorem, $H^0(\p^3,\cE(-N))=0$ for
for sufficiently large $N$, so that $H^0(\p^3,\cE(k)) = 0$ for all 
$k\leq-1$.

We also have that:
$$ 0 \to H^1(\p^3,\cE(k-1)) \to  H^1(\p^3,\cE(k)) ~ \forall k\leq-1 $$
hence $H^1(\p^3,\cE(-2)) = 0$ implies by induction that $H^1(\p^3,\cE(k)) = 0$
for all $k\leq-2$. Furthermore,
$$ H^2(\p^3,\cE(k-1)) \to  H^2(\p^3,\cE(k)) \to 0 ~ \forall k\geq-2 $$
forces $H^2(\p^3,\cE(k)) = 0$ for all $k\geq-2$ once $H^2(\p^3,\cE(-2)) = 0$.
\end{proof}

In this section we will concentrate on {\em framed admissible torsion-free sheaves},
that is a pair $(\cE,\phi)$ consisting of an admissible torsion-free sheaf such that
the restriction $\cE|_{\ell_\infty}$ is trivial plus a framing
$\phi:\cE|_{\ell_\infty}\stackrel{\sim}{\to}{\cal O}_{\ell_\infty}^{\oplus{\rm rk}\cE}$.
We will show that the moduli space of framed admissible torsion-free sheaves can
be parameterized by $\cpx$-stable complex ADHM data.
More about admissible sheaves in general can be found at \cite{J-p3}.

\subsection{From ADHM data to sheaves} \label{adhm2p3}
Let $(B_{kl},i_k,j_k)$ be a complex ADHM datum; combining constructions of
Donaldson \cite{D1} and Nakajima \cite{N2}, we define the monad:
\begin{equation} \label{monad.p3}
V\otimes\op3(-1) \stackrel{\alpha}{\longrightarrow}
\tW\otimes\op3 \stackrel{\beta}{\longrightarrow} V\otimes\op3(1)
\end{equation}
where the maps $\alpha$ and $\beta$ are given by:
\begin{equation} \label{alpha}
\alpha = \left( \begin{array}{c}
zB_{11} + wB_{21} + x \\ zB_{12} + wB_{22} + y \\ zj_1 + wj_2
\end{array} \right) \end{equation}
\begin{equation} \label{beta}
\beta = \left( \begin{array}{ccc}
-zB_{12} - wB_{22} - y \ \ & \ \ zB_{11} + wB_{21} + x \ \ & \ \ zi_1 + wi_2
\end{array} \right) \end{equation}

\begin{proposition} \label{l1}
$\beta\alpha=0$ if and only if $(B_{kl},i_k,j_k)$ satisfies the
complex ADHM equations (\ref{c1}-\ref{c3}).
\end{proposition}

The proof is a straightforward calculation left to the reader.
It follows from Proposition \ref{l1} that $\cE=\ker\beta/{\rm Im}\alpha$,
the first cohomology of the monad (\ref{monad.p3}), is a well-defined coherent
sheaf on $\p^3$.  We will now check that $\cE$ is the only nontrivial cohomology
of (\ref{monad.p3}). It is easy to see that $GL(V)$-equivalent complex ADHM
data will produce isomorphic cohomology sheaves.

\begin{proposition} \label{l2}
$\alpha_X$ is injective away from a subvariety of codimension 2.
\end{proposition}

In particular, $\alpha$ is injective as a sheaf map.

\begin{proof}
It easy to see that $\alpha$ is injective on the line 
$\ell_\infty=\{z=w=0\}$.
So consider a point $X=[x:y:z:w]\in\p^3\setminus\ell_{\infty}$, and
take $v\in V$ such that $\alpha_X(v)=0$, that is:
\begin{equation} \left\{ \begin{array}{l}
(zB_{11} + wB_{21})v = -xv \\ (zB_{12} + wB_{22})v = -yv \\
(zj_1 + wj_2)v = 0
\end{array} \right. \end{equation}
Thus $v$ is a common eigenvector of $zB_{11} + wB_{21}$ and
$zB_{12} + wB_{22}$, with eigenvalues $-x$ and $-y$, respectively. Hence,
for fixed $z,w\neq0$, we conclude that $x$ and $y$ may assume only finitely
many values. In other words, $\alpha_X$ may fail to be injective only at 
points
of the form $[x=x(z,w):y=y(z,w):z:w]$, so $\alpha_X$ is injective away from 
a codimension 2 subvariety (which does not intersect $\ell_\infty$).
\end{proof}

The following is the key result in the monad construction, and further
justifies our concept of $\cpx$-stability:

\begin{proposition} \label{l3}
$\beta$ is surjective if and only if $(B_{kl},i_k,j_k)$ is $\cpx$-stable.
\end{proposition}
\begin{proof}
Again, it is easy to see that $\beta$ is surjective on the line
$\ell_\infty=\{z=w=0\}$. So it is enough to show that the localization
of $\beta$ to all points $X=[x:y:z:w]\in\p^3\setminus\ell_{\infty}$
is surjective.

Equivalently, we argue that if $(B_{kl},i_k,j_k)$ is $\cpx$-stable,
then the dual map $\beta_X^*$ is injective for all 
$X\in\p^3\setminus\ell_{\infty}$.
Indeed, $\beta_X^*$ is not injective for some $[x:y:p_1:p_2]$, then there
is $v\in V$ such that:
\begin{equation} \left\{ \begin{array}{l}
\tilde{B}_1(p_1,p_2)^* v = \overline{x} v \\
\tilde{B}_2(p_1,p_2)^* v = -\overline{y} v \\
\tilde{\imath}(p_1,p_2)^* v = 0
\end{array} \right. \end{equation}
which, by duality, implies that
$\left( \tilde{B}_1(p_1,p_2),\tilde{B}_2(p_1,p_2),\tilde{\imath}(p_1,p_2) 
\right)$
is not stable. Thus $(B_{kl},i_k,j_k)$ is not $\cpx$-stable.

The converse statement is now clear: if $(B_{kl},i_k,j_k)$ is not 
$\cpx$-stable, then by duality $\beta_X^*$ is not injective for some
$[x:y:z:w]$, hence $\beta$ cannot be surjective.
\end{proof}

In order to further characterize the cohomology sheaf
$\cE=\ker\beta/{\rm Im}\alpha$, let $[H]$ denote the
generator of $H^\bullet(\p^3,\cpx)$, i.e. $[H]=c_1(\op3(1))$.

\begin{proposition} \label{l4}
The cohomology sheaf $\cE=\ker\beta/{\rm Im}\alpha$ is an admissible
torsion-free sheaf on $\p^3$, with ${\rm ch}(\cE)=r-c [H]^2$.
Moreover, $\cE|_{\ell_\infty}$ is trivial.
\end{proposition}
\begin{proof}
For the admissibility of $\cE$, set ${\cal K}=\ker\beta$. From the sequence
\begin{equation} \label{sqc1}
0 \to {\cal K} \to \tW\otimes\op3 \to V\otimes\op3(1) \to 0
\end{equation}
we obtain, after tensoring (\ref{sqc1}) by $\op3(k)$, the exact sequence of 
cohomology:
$$ V\otimes H^{p-1}(\p^3,\op3(k+1)) \to H^{p}(\p^3,{\cal K}(k)) \to 
\tW\otimes H^{p}(\p^3,\op3(k)) $$
which implies that $H^{p}(\p^3,{\cal K}(k))=0$ for $p\leq1$, $p+k\leq-1$ and
for $p\geq2$, $p+k\geq0$, since the two groups at the ends are zero in this 
range. Next, from the sequence:
\begin{equation} \label{sqc2}
0 \to V\otimes\op3(-1) \stackrel{\alpha}{\longrightarrow}
{\cal K} \longrightarrow \cE \to  0
\end{equation}
we obtain, after tensoring (\ref{sqc1}) by $\op3(k)$, the exact sequence
of cohomology:
$$ H^{p}(\p^3,{\cal K}(k)) \to H^{p}(\p^3,\cE(k)) \to V\otimes H^{p+1}(\p^3,\op3(k-1)) $$
and hence $H^{p}(\p^3,\cE(k))$ for $p\leq1$, $p+k\leq-1$ and for $p\geq2$, 
$p+k\geq0$, since the third group also vanishes in that range.

To compute ${\rm ch}(\cE)$, just notice that:
$$ {\rm ch}(\cE) = {\rm ch}(\tW\otimes\op3) - {\rm ch}(V\otimes\op3(-1))
- {\rm ch}(V\otimes\op3(1)) $$
since $\cE$ is the only non vanishing cohomology of the monad 
(\ref{monad.p3}). The triviality of $\cE|_{\ell_\infty}$ also follows
easily from the construction, see \cite{D1}.

It remains for us to show that $\cE$ is torsion-free. First,
notice that $\cal K$ is a locally-free sheaf, since $\beta_X$ is surjective
for all $X\in\p^3$ (see the proof of Proposition \ref{l3}). Moreover, as it
was pointed out in the proof of Proposition \ref{l2}, the $\alpha_X$ is
injective away from a subset of codimension 2 in $\p^3$. Applying
Proposition \ref{ap-prop2} in Appendix \ref{app-b} to sequence (\ref{sqc2}),
we conclude that $\cE$ must be torsion-free.
\end{proof}

We are finally in position to prove Proposition \ref{no} by looking at the
corresponding sheaves on $\p^3$.

\begin{proposition} \label{no.r=1}
There are no admissible torsion-free sheaves $\cE$ on $\p^3$ with
${\rm ch}(\cE) = 1 - c [H]^2$.
\end{proposition}

Indeed, if there were $\cpx$-stable solutions of (\ref{c1}-\ref{c3}) for 
$r=1$, the monad construction would produce admissible torsion-free
sheaves $\cE$ such that ${\rm ch}(\cE) = 1 - c [H]^2$. So the above
result implies Proposition \ref{no}.

\begin{proof}
First note that $\cE^{**}\simeq \op3$, since $\cE^{**}$ is reflexive of rank 
1 (hence locally-free) and $c_1(\cE^{**})=c_1(\cE)=0$. So $\cE$ is an
ideal sheaf fitting in the sequence
\begin{equation} \label{r1sqc}
0 \to \cE \to \op3 \to \cQ \to 0 ~~, ~~ \cQ = \op3/\cE ~ .
\end{equation}
Note that ${\rm ch}(\cQ) = c [H]^2$, and that $\cQ$ is the structure
sheaf of a 1-dimensional subscheme $s:\Sigma\hookrightarrow \p^3$
(i.e. $\cQ=s_*{\cal O}_{\Sigma}$).

After twisting the sequence (\ref{r1sqc}) by $\op3(k)$ and using 
admissibility, we get that $h^0(\p^3,\cQ(k))=0$ for all $k\leq-2$,
hence the $\Sigma={\rm supp}\cQ$ contains no 0-dimensional
components. Moreover, $h^1(\p^3,\cQ(k))=0$ for all $k\geq-2$; in
particular
$$ h^0(\Sigma,{\cal O}_{\Sigma})=h^0(\p^3,\cQ)=\chi(\p^3,\cQ)=2c~, $$
so that $\Sigma$ consists of $2c$ connected components 
$\Sigma_1,\cdots,\Sigma_{2c}$.

For each connected component $\Sigma_a$, we have that \linebreak
$\chi(\p^3,{\cal O}_{\Sigma_a})=\chi(\Sigma_a,{\cal O}_{\Sigma_a})=1$.
It follows that $\Sigma_a$ is a smooth $\p^1$, so $\cE$ must be the ideal 
sheaf of $2c$ lines in $\p^3$. But it is easy to check that the Chern character
of the ideal sheaf of $2c$ lines in $\p^3$ is given by $1-2c [H]^2+2c[H]^3$,
leading to a contradiction.
\end{proof}


\subsection{From sheaves to complex ADHM data}
The first step of the reverse construction is essentially
provided by Manin \cite{Ma}:

\begin{proposition} \label{l5}
Every admissible torsion-free sheaf $\cE$ on $\p^3$ can be obtained
as the cohomology of the monad
\begin{equation} \label{m1}
0\to  V\otimes\op3(-1) \stackrel{\alpha}{\longrightarrow}
\tW\otimes\op3 \stackrel{\beta}{\longrightarrow} V'\otimes\op3(1) \to 0 ~,
\end{equation}
where $V=H^1(\p^3,\cE\otimes\Omega^2_{\p^3}(1))$, 
$\tW=H^1(\p^3,\cE\otimes\Omega^1_{\p^3})$
and $V'=H^1(\cE(-1))$.
\end{proposition}
\begin{proof}
Manin proves the case $\cE$ being locally-free in \cite[p. 91]{Ma}, using
the Beilinson spectral sequence. However, the argument generalizes word by
word for $\cE$ being torsion-free; just note that the projection formula
$$ R^ip_{1*}\left(p_1^*\op3(k)\otimes p_2^*{\cal F}\right)=
\op3(k)\otimes H^i({\cal F}) $$
holds for every torsion-free sheaf $\cal F$, where $p_1$ and $p_2$ are the
natural projections of $\p^3\times\p^3$ onto the first and second factors.
\end{proof}

So let $\cE$ be an admissible torsion-free sheaf on $\p^3$ with
${\rm ch}(\cE) = r - c [H]^2$ and such that
$\cE|_{\ell_\infty}$ is trivial. It remains 
for us to show that the monad in Proposition \ref{l5} can be reduced
to a $\cpx$-stable solution of the complex ADHM equations
(\ref{c1}-\ref{c3}). A lengthy but straightforward cohomological 
calculation (see Appendix \ref{ap1}) shows that:
$$  h^1(\p^3,\cE\otimes\Omega^2_{\p^3}(1)) = h^1(\p^3,\cE(-1)) = c ~~,
~~ h^1(\p^3,\cE\otimes\Omega^1_{\p^3}) = c + 2r $$
and that there is a natural identification $V\simeq V'$.

Now, $\alpha\in{\rm Hom}(V,\tW)\h1$ and $\beta\in{\rm Hom}(\tW,V)\h1$, and 
we can express these maps in the following manner:
$$ \alpha =\alpha_1 x + \alpha_2 y + \alpha_3 z + \alpha_4 w ~~~{\rm and}~~~
\alpha =\beta_1 x + \beta_2 y + \beta_3 z + \beta_4 w $$
where, clearly, $\alpha_k\in{\rm Hom}(V,\tW)$ and
$\beta_k\in{\rm Hom}(\tW,V)$ for each $k=1,\dots,4$. The condition
$\beta\alpha=0$ then implies that:
$$ \beta_k\alpha_k = 0~,~~~k=1,\dots,4 $$
$$ \beta_k\alpha_l + \beta_l\alpha_k = 0~,~~~k,l=1,\dots,4~{\rm and}~k\neq l $$

Restricting (\ref{m1}) to the line at infinity $\ell_\infty = \{z=w=0\}$ we get:
$$ 0\to V\otimes\op3|_{\ell_\infty}(-1)
\stackrel{\alpha_\infty}{\longrightarrow}
\tW\otimes\op3|_{\ell_\infty}
\stackrel{\beta_\infty}{\longrightarrow}
V\otimes\op3|_{\ell_\infty}(1) \to 0 $$
where $\alpha_\infty=\alpha_1 x + \alpha_2 y$ and
$\beta_\infty=\beta_1 x + \beta_2 y$.
Setting ${\cal K}=\ker\beta$ we have:
$$ 0 \to V\otimes{\cal O}_{\ell_\infty}(-1) 
\stackrel{\alpha_\infty}{\longrightarrow}
{\cal K}|_{\ell_\infty} \longrightarrow \cE|_{\ell_\infty} \to 0 $$
from the associated long exact sequence of cohomology we conclude that
$H^1(\ell_\infty,{\cal K}|_{\ell_\infty})=0$ and
$H^0(\ell_\infty,{\cal K}|_{\ell_\infty})\simeq 
H^0(\ell_\infty,\cE|_{\ell_\infty})\simeq \cE_P$,
for some $P\in\ell_\infty$, since $H^p(\ell_\infty,{\cal 
O}_{\ell_\infty}(-1))=0$, for $p=1,2$,
and since $\cE|_{\ell_\infty}\simeq{\cal O}_{\ell_\infty}^{\oplus r}$. We 
set $W=H^0(\ell_\infty,{\cal K}|_{\ell_\infty})$; the choice of a basis for
$W$ corresponds to the choice of a trivialization for $\cE|_{\ell_\infty}$.

Similarly, from the sequence
$$ 0 \to {\cal K}|_{\ell_\infty} \longrightarrow \tW\otimes{\cal 
O}_{\ell_\infty}
\stackrel{\beta_\infty}{\longrightarrow}
V\otimes{\cal O}_{\ell_\infty}(1) \to 0 $$
we obtain:
\begin{equation} \label{m2}
0 \to W \longrightarrow \tW \stackrel{\beta_\infty}{\longrightarrow}
V\otimes H^0(\ell_\infty,{\cal O}_{\ell_\infty}(1)) \to 0
\end{equation}
since $H^0(\ell_\infty,{\cal O}_{\ell_\infty})\simeq\cpx$ and
$H^1(\ell_\infty,{\cal K}|_{\ell_\infty})=0$. Then using the identification
$H^0(\ell_\infty,{\cal O}_{\ell_\infty}(1))\simeq \cpx x \oplus \cpx y$
we can rewrite (\ref{m2}) in the following way:
\begin{equation} \label{wtw} 0 \to W \longrightarrow \tW
\stackrel{\left(\begin{array}{c} \beta_1 \\ \beta_2 
\end{array}\right)}{\longrightarrow}
V\oplus V \to 0 \end{equation}
so that $W = \ker\beta_1\cap\ker\beta_2$.

Applying the same argument to the dual monad:
$$ 0 \to V^*\otimes\op3|_{\ell_\infty}(-1) \stackrel{\beta_\infty^{\rm 
t}}{\longrightarrow}
\tW^*\otimes\op3|_{\ell_\infty}
\stackrel{\alpha_\infty^{\rm t}}{\longrightarrow}
V^*\otimes\op3|_{\ell_\infty}(1) \to 0 $$
we have the exact sequence:
$$ 0 \to H^0(\ell_\infty,\ker\{\alpha_\infty^{\rm t}\}) \longrightarrow 
\tW^*
\stackrel{\left(\begin{array}{c} \alpha_1 \\ \alpha_2 
\end{array}\right)}{\longrightarrow} V^*\oplus V^* $$
which implies that $(\alpha_1~~\alpha_2):V\oplus V\to \tW$ is injective. 
Moreover, the
sequence (\ref{wtw}) splits, and we can identify $\tW\simeq V\oplus V \oplus 
W$.

Furthermore, notice that
$$ \ker\beta_1/{\rm Im}\alpha_1\simeq \cE_{[1,0,0,0]}\simeq W \simeq 
\ker\beta_1\cap\ker\beta_2 ~.$$
Thus ${\rm Im}\alpha_1\cap\ker\beta_2=0$, so that 
$\beta_1\alpha_2=-\beta_2\alpha_1:V\to V$ are
isomorphisms.

Therefore we have:
$$ \alpha_1 = \left( \begin{array}{c}
\id_V \\ 0 \\ 0 \end{array} \right) ~,~
\alpha_2 = \left( \begin{array}{c}
0 \\ \id_V \\ 0 \end{array} \right) ~~~
\begin{array}{c} \beta_1 = \left( \begin{array}{ccc}
0 ~~ & ~~ \id_V ~~ & ~~ 0
\end{array} \right) \\ \\
\beta_2 = \left( \begin{array}{ccc}
-\id_V ~~ & ~~ 0 ~~ & ~~ 0
\end{array} \right) \end{array} $$
and the condition $\beta\alpha=0$ implies that:
$$ \alpha_3 = \left( \begin{array}{c}
B_{11} \\ B_{12} \\ j_1 \end{array} \right) ~,~
\alpha_4 = \left( \begin{array}{c}
B_{21} \\ B_{22} \\ j_2 \end{array} \right) ~,~
\begin{array}{c} \beta_3 = \left( \begin{array}{ccc}
-B_{12} ~~ & ~~ B_{11} ~~ & ~~ i_1
\end{array} \right) \\ \\
\beta_4 = \left( \begin{array}{ccc}
-B_{22} ~~ & ~~ B_{21} ~~ & ~~ i_2
\end{array} \right) \end{array} $$
with $(B_{kl},i_k,j_k)$ being a complex ADHM datum satisfying the complex
ADHM equations (\ref{c1}-\ref{c3}). The surjectivity of $\beta$ implies the
$\cpx$-stability of $(B_{kl},i_k,j_k)$, by Proposition \ref{l3}. Summing up, 
we have proved:

\begin{theorem}\label{MT2}
There is a 1-1 correspondence between the following objects:
\begin{itemize}
\item framed admissible torsion-free sheaves on $\p^3$, and 
\item $\cpx$-stable solutions of the complex ADHM equations.
\end{itemize}
In particular, the moduli space of framed admissible torsion-free sheaves $\cE$
on $\p^3$ with ${\rm ch}(\cE)=r- c [H]^2$ is a smooth complex manifold of
dimension $4rc$, non-empty for $r\geq2$.
\end{theorem}


\subsection{Locally-free admissible sheaves and complex instantons} 
\label{4.3}

We will now describe necessary and sufficient conditions that guarantee
that the cohomology sheaf of the monad (\ref{monad.p3}) is locally-free.
Recall that $\alpha_X$ denotes the localization of the map $\alpha$
to a point $X\in\p^3$.

\begin{proposition} \label{l2'}
$\alpha_X$ is injective for all $X\in\p^3$ if and only if
$(B_{kl},i_k,j_k)$ is $\cpx$-costable.
\end{proposition}
\begin{proof}
Since $\alpha$ is injective on the line
$\ell_\infty=\{z=w=0\}$, it is enough to show that $\alpha_X$ is
injective for all $X=[x:y:z:w]\in\p^3\setminus\ell_{\infty}$.

Indeed, take $v\in V$ such that $\alpha_X(v)=0$, hence:
$$ \left\{ \begin{array}{l}
\tilde{B}_1v = -xv \\ \tilde{B}_2v = -yv \\
\tilde{\jmath}v = 0
\end{array} \right. $$
But $(\tilde{B}_k,\tilde{\jmath},\tilde{\jmath})$ is costable for all
$(z,w)\in\cpx^2\setminus\{0\}$, therefore $v=0$.

Conversely, if $(B_{kl},i_k,j_k)$ is not $\cpx$-costable, there are
$(\lambda,\mu)\in\cpx^2\setminus\{0\}$ and a proper subspace $S\subset V$ 
such that $ [ \tilde{B}_1,\tilde{B}_2 ] |_S = 0$ and 
$S\subset\ker\tilde{\jmath}$.
Therefore $\alpha_{[x:y:\lambda:\mu]}(v)=0$ for all $v\in S$.
\end{proof}

Thus if $\vB=(B_{kl},i_k,j_k)$ is $\cpx$-regular, then $\alpha_X$ is 
injective and $\beta_X$ is surjective for all $X\in\p^3$, so that the quotient
$\ker\beta_X / {\rm Im} \alpha_X$ is a vector space of dimension $r$ for all
$X\in\p^3$. We conclude that:

\begin{corollary}
The cohomology sheaf $\cE$ is locally-free if and only if
$(B_{kl},i_k,j_k)$ is $\cpx$-regular.
\end{corollary}

\begin{remark} \rm
As it was pointed out in Remark \ref{st.not.reg}, there are solutions of the 
complex ADHM equations which are $\cpx$-regular but not $\cpx$-stable.
Therefore, {\em there exist admissible torsion-free sheaves which are not locally-free}.
The basic example is the cohomology $\cE$ of the monad:
$$ \op3(-1) \stackrel{\alpha}{\rightarrow} \op3^{\oplus4}
\stackrel{\beta}{\rightarrow} \op3(1) $$
$$ \alpha = \left(\begin{array}{c} x \\ y \\ 0 \\ 0 \end{array}\right) 
~~{\rm and}~~
\beta= (-y ~~ x ~~ z ~~ w) $$
It is easy to see that $\beta$ is surjective for all $(x,y,z,w)\in\p^3$, 
while $\alpha$ is injective provided $x,y\neq0$. It then follows from applying
Proposition \ref{ap-prop2} to sequence (\ref{sqc2}) that $\cE$ is 
torsion-free, but not locally-free. In particular, the singularity set of $\cE$
(i.e. the support of $\cE^{**}/\cE)$ consists of the line $\{x=y=0\}\subset\p^3$.
\end{remark}

Rank 2 locally-free sheaves $\cE$ with $c_1(\cE)=0$ and $H^0(\p3,E(-1))=H^1(\p^3,\cE(-2))=0$
are also known in the literature as {\em mathematical} (or {\em complex}) {\em instanton bundles}
(see \cite{OSS,WW} and also \cite{CTT} for more recent references and a brief
survey of the subject); it is easy to see, via Serre duality, that these are admissible.
They correspond, via the Penrose transform, to holomorphic vector bundles with
$SL(2,\cpx)$ anti-self-dual connections on $\mathbb{M}$, the complexified
compactified Minkowski space-time (see \cite{WW}; recall also that $\mathbb{M}$
is just the Grassmaniann of lines in $\p^3$). The integer $c=c_2(\cE)$ is also
called the {\em charge} of the complex instanton bundle $\cE$. 

For rank $r>2$, a complex instanton bundle is an admissible locally-free sheaf
and a framed complex instanton bundle is a framed admissible locally-free sheaf.
Clearly, the moduli space of equivalence classes of framed complex instanton
bundles fibers over the moduli space of equivalence classes of complex instanton
bundles, with fibers given by $PGL(r,\cpx)$, the set of all possible framings.
It follows that the moduli space of framed complex instanton bundles of rank
$r\geq2$ and charge $c\geq1$ is exactly ${\cal M}_{\cpx}^{\rm reg}(r,c)$, the
open subset of $\cmc$ consisting of the orbits of $\cpx$-regular solutions of the 
complex ADHM equations (\ref{c1}-\ref{c3}). 

Determining the irreducibility and smoothness of the moduli space of rank 2 complex
instanton bundles with charge $c$ is a long standing question, see \cite{CTT} for
a recent short survey of this topic. As a
special case of Theorem \ref{MT2}, we obtain a strong result along these lines:

\begin{corollary}
The moduli space of framed complex instanton bundles of rank $r\geq2$
and charge $c$ is a nonempty, smooth complex manifold of dimension $4rc$.
\end{corollary}

We remark that framed complex instanton bundles are always $\mu$-semistable,
see \cite[p. 210]{OSS}. It would be interesting to compare the admissibility and
semistability condition, and determine under what necessary and sufficient conditions
admissible torsion-free sheaves are $\mu$-semistable, and vice-versa; a few results
along these lines have been obtained by the first author in \cite{J-p3}.  


\subsection{Reflexive admissible sheaves}

Reflexive sheaves on $\p^3$ have been extensively studied in a series of 
papers by Hartshorne \cite{Ha}, among other authors. In particular, it was
show that a rank 2 reflexive sheaf on $\p^3$ is locally-free if and only
if $c_3({\cal F})=0$. Therefore, we conclude:

\begin{proposition} ({\bf Hartshorne \cite{Ha}})
There are no rank 2 admissible sheaves on $\p^3$ which are reflexive but not 
locally-free.
\end{proposition}

The situation for higher rank is quite different, though, and it is easy to 
construct a rank 3 admissible sheaf which is reflexive but not locally-free.
Setting $r=3$ and $c=1$, consider the monad:
$$ \op3(-1) \stackrel{\alpha}{\rightarrow} \op3^{\oplus5}
\stackrel{\beta}{\rightarrow} \op3(1) $$
$$ \alpha = \left(\begin{array}{c} x \\ y \\ 0 \\ 0 \\ z \end{array}\right) 
~~{\rm and}~~
\beta= (-y ~~ x ~~ z ~~ w ~~ 0) $$
Again, it is easy to see that $\beta$ is surjective for all $(x,y,z,w)\in\p^3$,
while $\alpha$ is injective provided $x,y,z\neq0$. It then follows from applying
Proposition \ref{ap-prop2} to sequence (\ref{sqc2}) that $\cE$ is reflexive,
but not locally-free; its singularity set is just the point $[0:0:0:1]\in\p^3$.

\begin{proposition}
The cohomology sheaf $\cE$ is reflexive if and only if
$\vB$ is $\cpx$-semiregular.
\end{proposition}

It then follows from Hartshorne's result that every $\cpx$-semiregular solution of
the complex ADHM equations for $r=2$ is $\cpx$-regular, as we claimed in
Remark \ref{st.not.reg}; the example of a properly reflexive admissible sheaf
corresponds to the properly $\cpx$-semiregular solution of the complex ADHM
equations for $r=3$ given in Remark \ref{st.not.reg}.

\begin{proof}
If $\vB$ is $\cpx$-semiregular, then $\vB_p$ is costable for all but 
finitely many $p\in\p^1$ (see the observation following Proposition
\ref{ss}). This means that $\alpha_X$ is injective for all but finitely
many $X\in\p^3$. Thus, by Proposition \ref{ap-prop2}, the cohomology
sheaf $\cE$ is reflexive.

Conversely, if $\cE$ is reflexive then $\alpha_X$ is injective for all but 
finitely many $X\in\p^3$. It follows that $\vB_p$ must be costable for all
but finitely many $p\in\p^1$, and $\vB$ is $\cpx$-semiregular.
\end{proof}


\section{Quantum instantons} \label{qinst}

In this section we will adapt the ADHM construction of instantons to obtain
{\em complex quantum instantons} ({\em cf.} \cite{FJ}) from  $\cpx$-regular
solutions of the complex ADHM equations. For the convenience of the reader,
let us first recall the essential definitions from our previous paper.

\subsection{Quantum Minkowski space-time}  \label{qinst1}

In \cite{FJ}, we defined the {\em quantum compactified, complexified 
Minkowski space} $\mpq$ as the associative graded $\cpx$-algebra generated by
$z_{11'},z_{12'},z_{21'},z_{22'},D,D'$ satisfying the relations
(\ref{relations1}) to (\ref{relations5}) below ($p=q^{\pm 1}$ are formal 
parameters):

\begin{equation} \label{relations1}
\begin{array}{lcr}
z_{11'}z_{12'}=z_{12'}z_{11'} & \ \ \ & z_{11'}z_{21'}=z_{21'}z_{11'} \\
z_{12'}z_{22'}=z_{22'}z_{12'} & \ \ \ & z_{21'}z_{22'}=z_{22'}z_{21'} \\
& z_{12'}z_{21'}=z_{21'}z_{12'} &
\end{array} \end{equation}
\begin{equation} \label{relations2}
q^{-1} (z_{11'}z_{22'}-z_{12'}z_{21'})=q (z_{22'}z_{11'}-z_{12'}z_{21'})
\end{equation}
\begin{equation} \label{relations3}
\begin{array}{lcr}
Dz_{11'}=pq^{-1} z_{11'}D & \ \ \ & D'z_{11'}=p^{-1}q^{-1} z_{11'}D' \\
Dz_{12'}=pq^{-1} z_{12'}D & \ \ \ & D'z_{12'}=p^{-1}q z_{12'}D'      \\
Dz_{21'}=pq z_{21'}D      & \ \ \ & D'z_{21'}=p^{-1}q^{-1} z_{21'}D' \\
Dz_{22'}=pq z_{22'}D      & \ \ \ & D'z_{22'}=p^{-1}q z_{22'}D'      \\
\end{array} \end{equation}
\begin{equation} \label{relations4}
p^{-1} DD'=p D'D
\end{equation}
\begin{equation} \label{relations5}
q^{-1} (z_{11'}z_{22'}-z_{12'}z_{21'}) = p^{-1} DD'
\end{equation}

Localizations of $\mpq$ with respect to $D$ and $D'$ lead to two ``affine 
patches'' $\miq$ and $\mjq$, respectively. More precisely, the new generators:
$$ x_{rs'} = \frac{z_{rs'}}{D} ~~~ {\rm and} ~~~ y_{rs'} = 
\frac{z_{rs'}}{D'} $$
satisfy the relations:
\begin{equation} \label{xcommut1}
x_{11'}x_{12'} = x_{12'}x_{11'}, ~~~ x_{21'}x_{22'} = x_{22'}x_{21'},
\end{equation}
\begin{equation} \label{xcommut3}
[x_{11'},x_{22'}] + [x_{21'},x_{12'}] = 0
\end{equation}
\begin{equation} \label{xcommut2} \begin{array} {c}
x_{11'}x_{21'} = q^{-2}x_{21'}x_{11'}, \ \ \ x_{12'}x_{22'} =
q^{-2}x_{22'}x_{12'}, \\
x_{21'}x_{12'} = q^2 x_{12'} x_{21'}
\end{array} \end{equation}
and
\begin{equation} \label{ycommut1}
y_{11'}y_{21'} = y_{21'}y_{11'}, ~~~ y_{12'}y_{22'} = y_{22'}y_{12'},
\end{equation}
\begin{equation} \label{ycommut3}
\left[ y_{11'},y_{22'} \right] + [y_{12'},y_{21'}] = 0
\end{equation}
\begin{equation} \label{ycommut2} \begin{array} {c}
y_{11'}y_{12'} = q^{-2}y_{12'}y_{11'}, \ \ \ y_{21'}y_{22'} = 
q^{-2}y_{22'}y_{21'}, \\
y_{12'}y_{21'} = q^2 y_{21'} y_{12'}
\end{array} \end{equation}
$$ {\rm i.e.} ~~ \miq = \mpq[D^{-1}] = 
\cpx[x_{11'},x_{12'},x_{21'},x_{22'}]/(\ref{xcommut1}-\ref{xcommut2})~. $$
$$ {\rm and} ~~ \mjq = \mpq[D'^{-1}] = 
\cpx[y_{11'},y_{12'},y_{21'},y_{22'}]/(\ref{ycommut1}-\ref{ycommut2})~. $$

These two algebras can be made isomorphic after inverting the determinants
$\det(x)=x_{11'}x_{22'}-x_{12'}x_{21'}$ and $\det(y)=y_{11'}y_{22'}-y_{21'}y_{12'}$.
Note that $\det(x)=\det(y)^{-1}=\frac{D'}{D}$. One has:
\begin{equation} \label{nc.glueing}
\miq[\det(x)^{-1}] \stackrel{\sim}{\longrightarrow} \mjq[\det(y)^{-1}]
\end{equation}
$$ {\rm where} ~~ y_{kl'} = \frac{x_{kl'}}{\det(x)} ~.$$
We denote $\mijq$ the isomorphic algebras in (\ref{nc.glueing}).

Let us now focus on $\miq$; all observations below will also apply to $\mjq$.

It follows immediately from the commutation relations 
(\ref{xcommut1}-\ref{xcommut2}) that any element of $\miq$ can
be presented as a sum of monomials of the form:
\begin{equation} \label{basis}
x_{11'}^{n_{11'}} x_{12'}^{n_{12'}} x_{21'}^{n_{21'}} x_{22'}^{n_{22'}} ~~ , 
~~ n_{ij'}\geq0
\end{equation}
Moreover, it is easy to see directly from (\ref{xcommut1}-\ref{xcommut2}), 
and it is also proven in \cite[Theorem 21]{FJ}, that these monomials are linearly 
independent and therefore form a basis of $\miq$. An element $f\in\miq$ has degree
$d$ if it is a sum of monomials (\ref{basis}) with $n_{11'} + n_{12'} + n_{21'} + n_{22'} = d$.

A concise form of the commutation relations (\ref{xcommut1}-\ref{xcommut2}) 
can also be expressed in terms of an $R$-matrix:
\begin{equation} \label{rmatrixI}
R^{\rm I}_{12} = \left( \begin{array}{cccc}
p^{-1} & 0 & 0 & 0 \\ 0 & q^{-1} & p^{-1}-q & 0 \\
0 & p^{-1}-q^{-1} & q & 0 \\ 0 & 0 & 0 & p^{-1}
\end{array} \right), \ \ p=q^{\pm1}
\end{equation}
Setting
$$ X = \left( \begin{array}{cc}
x_{11'} & x_{12'} \\ x_{21'} & x_{22'}
\end{array} \right) $$
and defining $X_1=X\otimes\id$ and $X_2=\id\otimes X$, the relations
(\ref{xcommut1}-\ref{xcommut2}) become equivalent to the identity:
$$ R^{\rm I}_{12} X_1X_2 =  X_2X_1 R^{\rm I}_{12} $$

We define the module of 1-forms over $\miq$, denoted by $\Omega^1_{\miq}$, 
as the $\miq$-bimodule
generated by:
$$ dX = \left( \begin{array}{cc}
dx_{11'} & dx_{12'} \\ dx_{21'} & dx_{22'}
\end{array} \right) $$
satisfying the following relations (written in matrix form):
\begin{equation} \label{xnc1forms}
R^{\rm I}_{12} X_1 dX_2 = dX_2  X_1 (R^{\rm I}_{21})^{-1}
\end{equation}
where $dX_2 = \id\otimes dX$ and
$R^{\rm I}_{21}=Q_1^{-1}Q_2R_{21}Q_1^{-1}Q_{2}$

Similarly, the module of 2-forms $\Omega^2_{\miq}$ is the $\miq$-bimodule 
generated by $dx_{rs'}\wedge dx_{kl'}$ satisfying the relations (written in
matrix form):
\begin{equation} \label{xnc2forms}
R^{\rm I}_{12} dX_1 \wedge dX_2 = - dX_2 \wedge dX_1 (R^{\rm I}_{21})^{-1}
\end{equation}
where $dX_1 = dX \otimes \id$. The module $\Omega^2_{\miq}$ splits as the 
sum of the submodules:
\begin{equation} \label{nc-sd2f}
\Omega^{2,+}_{\miq} = \miq \langle dx_{11'}\wedge dx_{12'},
dx_{21'}\wedge dx_{22'}, dx_{11'}\wedge dx_{22'} - dx_{12'}\wedge dx_{21'}
\rangle \end{equation}
\begin{equation} \label{nc-asd2f}
\Omega^{2,-}_{\miq} = \miq \langle dx_{11'}\wedge dx_{21'},
dx_{12'}\wedge dx_{22'}, dx_{11'}\wedge dx_{22'} + dx_{12'}\wedge dx_{21'}
\rangle \end{equation}
which can be regarded as the modules of self-dual and anti-self-dual 
2-forms, respectively.

Finally, the action of the de Rham operator $d: \miq\to\Omega^1_{\miq}$ is 
given on the generators as $x_{rs'}\mapsto dx_{rs'}$, and it is then extended
to the whole $\miq$ by $\cpx$-linearity and the Leibnitz rule:
\begin{equation} \label{l}
d(fg) = g df + f dg
\end{equation}
where $f,g \in \miq$. One also defines the de Rham operator $d : 
\Omega^1_{\miq}\to\Omega^2_{\miq}$ on the
generators as $f dx_{rs'}\mapsto df \wedge dx_{rs'}$, also extending it by 
$\cpx$-linearity and the Leibnitz rule (\ref{l}). Relations (\ref{xnc1forms})
and (\ref{xnc2forms}) imply that $d^2=0$.

Now let $E$ be a right $\miq$-module. A {\em connection} on $E$ is a 
$\cpx$-linear map:
$$ \nabla: E \to E\otimes_{\miq}\Omega^1_{\miq} $$
satisfying the Leibnitz rule:
$$ \nabla(\sigma f)=\sigma\otimes df + \nabla(\sigma) f $$
where $f\in\miq$ and $\sigma\in E$. The connection $\nabla$ also acts on 
1-forms, being defined as the $\cpx$-linear map:
$$ \nabla: E\otimes_{\miq}\Omega^1_{\miq} \to
E\otimes_{\miq}\Omega^2_{\miq} $$
satisfying:
$$ \nabla(\sigma\otimes\omega)=\sigma\otimes d\omega +
\nabla\sigma \wedge \omega $$
where $\omega\in\Omega^1_{\miq}$.

Moreover, two connections $\nabla$ and $\nabla'$ are said to be gauge 
equivalent if there is $g\in{\rm Aut}_{\miq}(E)$ such that
$\nabla = g^{-1} \nabla' g$.

The {\em curvature} $F_\nabla$ is defined by the composition:
$$ E \stackrel{\nabla}{\longrightarrow} E\otimes_{\miq}\Omega^1_{\miq}
\stackrel{\nabla}{\longrightarrow} E\otimes_{\miq} \Omega^2_{\miq} $$
and it is easy to check that $F_\nabla$ is actually right $\miq$-linear. 
Therefore, $F_\nabla$
can be regarded as an element of ${\rm 
End}_{\miq}(E)\otimes_{\miq}\Omega^2_{\miq}$.
Furthermore, if $\nabla$ and $\nabla'$ are gauge equivalent, then there is 
$g\in{\rm Aut}_{\miq}(E)$
such that $F_\nabla = g^{-1} F_{\nabla'} g$. A connection $\nabla$ is said 
to be anti-self-dual if
$F_\nabla\in {\rm End}_{\miq}(E)\otimes_{\miq}\Omega^{2,-}_{\miq}$.

\begin{definition}
A {\em complex quantum instanton} over $\mpq$ consists of the following data:
\begin{enumerate}
\item finitely generated free right $\miq$- and $\mjq$-modules $E_{\rm I}$ and $E_{\rm J}$;
\item anti-self-dual connections $\nabla_{\rm I}$ and $\nabla_{\rm J}$ on 
$E_{\rm I}$ and $E_{\rm J}$, respectively;
\item an isomorphism $\Gamma : E_{\rm I}[\det(x)^{-1}] \to E_{\rm J}[\det(y)^{-1}]$ 
satisfying $\nabla_{\rm J}\Gamma=\Gamma\nabla_{\rm I}$.
\end{enumerate}\end{definition}

The attentive reader will notice that in our previous paper \cite{FJ} we
defined a quantum instanton as a pair consisting of a projective module
and an anti-self-dual connection. As we will see below, the modules
produced via ADHM construction are actually free, and we will use the above
definition in this paper. The construction of projective modules which are
not free is an interesting direction for future research.


\subsection{Construction of complex quantum instantons} \label{qinst2}
We will use a variation of the celebrated ADHM construction of instantons 
\cite{ADHM} to construct complex quantum instantons from $\cpx$-stable
solutions of the complex ADHM equations (\ref{c1}-\ref{c3}).

To begin, let $(B_{kl},i_k,j_k)\in\vB$ be a complex ADHM datum, and consider 
the maps:
$$ \xymatrix{ V\otimes\miq \ar[r]^{\alpha_1}_{\alpha_2} &
\tW\otimes\miq \ar[r]^{\beta_1}_{\beta_2} & V\otimes\miq } $$
defined as follows:
$$ \begin{array}{lcr}
\alpha_1=\left( \begin{array}{c}
B_{11}\otimes\id - \id\otimes x_{11'} \\ B_{12}\otimes\id - \id\otimes 
x_{12'} \\ j_1\otimes\id
\end{array} \right) & {\rm and} &
\alpha_2=\left( \begin{array}{c}
B_{21}\otimes\id - \id\otimes x_{21'} \\ B_{22}\otimes\id - \id\otimes 
x_{22'} \\ j_2\otimes\id
\end{array} \right) \end{array} $$
$$ \beta_1=\left( \begin{array}{lcr}
-B_{12}\otimes\id + \id\otimes x_{12'} \ \ & \ \ B_{11}\otimes\id - 
\id\otimes x_{11'} \ \ & \ \ i_1\otimes\id
\end{array} \right) $$
$$ \beta_2=\left( \begin{array}{lcr}
-B_{22}\otimes\id + \id\otimes x_{22'} \ \ & \ \ B_{21}\otimes\id - 
\id\otimes x_{21'} \ \ & \ \ i_2\otimes\id
\end{array} \right) $$

\begin{proposition} \label{ids}
$(B_{kl},i_k,j_k)$ satisfies the complex ADHM equations (\ref{c1}-\ref{c3})
if and only if the following identities hold:
\begin{equation} \label{c12'} \beta_1\alpha_1 = \beta_2\alpha_2 = 0 
\end{equation}
\begin{equation} \label{c3'} \beta_2\alpha_1 + \beta_1\alpha_2 = 0 
\end{equation}
\end{proposition}
\begin{proof}
It is easy to check that:
$$ \beta_1\alpha_1 = 
([B_{11},B_{12}]+i_1j_1)\otimes\id+\id\otimes[x_{11'},x_{12'}] $$
$$ \beta_2\alpha_2 = 
([B_{21},B_{22}]+i_2j_2)\otimes\id+\id\otimes[x_{21'},x_{22'}] $$
$$ \beta_2\alpha_1 + \beta_1\alpha_2 =
([B_{11},B_{22}]+[B_{21},B_{12}]+i_1j_2+i_2j_1)\otimes\id+\id\otimes([x_{11'},x_{22'}]+[x_{21'},x_{12'}])  $$
so that the statement follows easily from the commutation relations 
(\ref{xcommut1}) and (\ref{xcommut3}).
\end{proof}

For points $P=[p_1:p_2]$ and $Q=[q_1:q_2]$ in $\p^1$, we consider maps:
$$ \beta_P = p_1\beta_1 + p_2\beta_2 : \tW\otimes\miq \to V\otimes\miq $$
$$ \alpha_Q = q_1\alpha_1 + q_2\alpha_2 : V\otimes\miq \to \tW\otimes\miq $$
It is easy to see that if $\vec{B}$ satisfies (\ref{c12'}-\ref{c3'}), then:
\begin{equation} \label{bpaq}
\beta_P\alpha_Q = (p_1q_2 - p_2q_1) \beta_1\alpha_2
\end{equation}

\begin{proposition} \label{inj-surj}
Assume that $\vec{B}$ satisfies the complex ADHM equations 
(\ref{c12'}-\ref{c3'}).
\begin{enumerate}
\item $\beta_P\alpha_Q$ is injective for all $P\neq Q\in\p^1$.
\item $\beta_P$ is surjective if and only if $\vec{B}_P$ is stable.
\end{enumerate}
In particular, $\beta_P$ is surjective $\forall~P\in\p^1$ if and only if 
$\vec{B}$ is $\cpx$-stable.
\end{proposition}

It also follows easily that $\alpha_Q$ is injective $\forall~Q\in\p^1$, so 
that ${\rm Im}~\alpha_Q$ is a free submodule of $\tW\otimes\miq$, of rank $c$.
Furthermore,  $\ker\beta_P\cap{\rm Im}~\alpha_Q=\{0\}$ for all $P\neq 
Q\in\p^1$.

\begin{proof}
For the first statement, it is enough to show that $\beta_1\alpha_2$ is 
injective, by (\ref{bpaq}). Note that
any $\nu\in V\otimes\miq$ can be presented as a sum 
$\nu=\nu_d+\nu_{d-1}+\cdots+\nu_0$ where
$$ \nu_d = \sum_k v_k \otimes f_k ~~ , ~~ v_k\in V, f_k \in\miq ~ {\rm 
has~degree} ~ d ~ ,$$
so that $\beta_1\alpha_2 (\nu) = (x_{11'}x_{22'}-x_{12'}x_{21'})\nu_d +$(terms of lower degree).

We argue that the endomorphism of $\miq$ given by the multiplication by
$(x_{11'}x_{22'}-x_{12'}x_{21'})$ is injective. In fact, we can present an 
element $f\in\miq$ in the basis of monomials (\ref{basis}). Let us choose the
lexicographic order of such basis, i.e.
$(n_{11'},n_{12'},n_{21'},n_{22'})>(m_{11'},m_{12'},m_{21'},m_{22'})$
if $n_{11'}>m_{11'}$, or $n_{11'}=m_{11'}$ and $n_{12'}>m_{12'}$, or
$n_{11'}=m_{11'}$ and $n_{12'}=m_{12'}$ and $n_{21'}>m_{21'}$, or
$n_{11'}=m_{11'}$ and $n_{12'}=m_{12'}$ and $n_{21'}=m_{21'}$ and 
$n_{22'}>m_{22'}$. Then the commutation relations (\ref{xcommut1}-\ref{xcommut3})
imply that, in this basis, the multiplication by $x_{ij'}$ increases the exponent 
$n_{ij'}$ by $1$ and multiplies its coefficient by a power of $q$. Thus the
multiplication by $(x_{11'}x_{22'}-x_{12'}x_{21'})$ of a polynomial with a
nonzero leading monomial of the form (\ref{basis}) will yield a polynomial
with a nonzero leading monomial of the form
$x_{11'}^{n_{11'}+1} x_{12'}^{n_{12'}} x_{21'}^{n_{21'}} x_{22'}^{n_{22'}+1}$.
Thus indeed, multiplication by $(x_{11'}x_{22'}-x_{12'}x_{21'})$ induces an 
injective endomorphism of $\miq$, as desired.

Now, if $\beta_1\alpha_2 (\nu) = 0$, then $(x_{11'}x_{22'}-x_{12'}x_{21'})\nu_d = 0$,
which in turn implies that $\nu_d = 0$; by induction on $d$, we conclude that $\nu=0$.

For the second statement, consider the following polynomial algebras for 
each $[p_1:p_2]\in\p^1$:
$$ \chi_P = \cpx[p_1x_{11'} + p_2x_{21'} , p_1x_{12'} + p_2x_{22'}] ~~ {\rm 
and} ~~ \chi_{\overline{P}} = \cpx[p_2x_{11'} - p_1x_{21'} , p_2x_{12'} - p_1x_{22'}] $$
as commutative subalgebras of $\miq$; clearly, $\miq=  \chi_P \otimes \chi_{\overline{P}}/\sim$,
where by "$\sim$'' we understand the commutation relations between the generators of $\chi_P$
and those of  $\chi_{\overline{P}}$, which can be deduced from (\ref{xcommut1}-\ref{xcommut3}).

As in \cite[Proposition 10]{FJ}, we see that, for each $P\in\p^1$, $\beta_P$
restricts to a map $\beta_P|_{\chi_P}:\tW\otimes\chi_P\to V\otimes\chi_P$, and that
$\beta_P=\beta_P|_{\chi_P}\otimes\id_{\chi_{\overline{P}}}$. To complete the proof,
recall from \cite[Lemma 2.7]{N2} that $\beta_P|_{\chi_P}$ is surjective if and only if 
$\vec{B}_P$ is stable.
\end{proof}

Now we consider the maps $\vec{\alpha}:(V\oplus V)\otimes\miq\to\tW\otimes\miq$
and $\vec{\beta}:\tW\otimes\miq\to (V\oplus V)\otimes\miq$ given by:
$$ \vec{\alpha} = \left( \begin{array}{lr} \alpha_1 \ \ & \ \ \alpha_2 
\end{array} \right)
~~~ {\rm and} ~~~ \vec{\beta} = \left( \begin{array}{c} -\beta_2 \\ \beta_1 
\end{array} \right) $$
If $\vec{B}$ satisfies the complex ADHM equations, the identities in 
Proposition \ref{ids} imply that $\Xi=\vec{\beta}\vec{\alpha}= \beta_1\alpha_2 \id_{\cpx^2}$.
It follows that $\Xi$ is injective; in particular, $\vec{\alpha}$ is 
injective and ${\rm Im}~\vec{\alpha}={\rm Im}~\alpha_1\oplus{\rm Im}~\alpha_2$ is a
free submodule of $\tW\otimes\miq$, of rank $2c$. Moreover,
$\ker\vec{\beta}\cap{\rm Im}~\vec{\alpha}=\{0\}$.

\begin{proposition} \label{proj}
If $\vec{B}$ is a $\cpx$-regular solution of the complex ADHM equations, then
the map $\beta_1\alpha_2:V\otimes\miq\to V\otimes\miq$ is an isomorphism; if
$\beta_1\alpha_2$ is an isomorphism, then $\vec{B}$ is $\cpx$-stable. 
\end{proposition}

\begin{proof}
We already know that $\beta_1\alpha_2$ is injective. If $\vec{B}$ is $\cpx$-regular, then 
$\beta_1\alpha_2$ is also surjective in the classical case $q=1$. Since
$\dim{\rm coker}~\beta_1\alpha_2$ cannot jump for generic value of the parameter
$q$, we obtain the first statement. 

Now if $\beta_1\alpha_2$ is an isomorphism then, according to (\ref{bpaq}),  $\beta_P$ is
surjective for all $P\in\p^1$; thus $\vec{B}$ is a $\cpx$-stable by Proposition \ref{inj-surj}.
\end{proof}

In particular, if $\vec{B}$ is a $\cpx$-regular solution of (\ref{c1}-\ref{c3}),
then $\Xi=\vec{\beta}\vec{\alpha}$ is an isomorphism, and we define the map:
$$ P: \tW\otimes\miq \to \tW\otimes\miq $$
$$ P = \id_{\tW\otimes\miq} - \vec{\alpha} (\Xi)^{-1} \vec{\beta} $$
and notice that $P^2=P$, i.e. $P$ is a projection. Note that:
$$ {\rm Im}(P) = \ker\vec{\beta} = \ker\beta_1~\cup~\ker\beta_2 ~ . $$
The right $\miq$-module $E={\rm Im}(P)$ is finitely generated and 
stably-free, since $\ker P={\rm Im}~\vec{\alpha}$ is free ($\vec{\alpha}$ is
injective) and $E\oplus\ker P=\tW\otimes\miq$. Furthermore, $\miq$ is Noetherian
and every stably-free module over a noetherian ring is free. Thus we conclude
that $E$ is a free $\miq$-module.

A connection $\nabla$ on $E$ can be easily defined via the projection formula, as usual:
$$ \xymatrix{ \nabla: E \ar[r]^{\iota} & \tW\otimes\miq \ar[r]^{\id\otimes d} &
\tW\otimes\Omega_{\miq}^1 \ar[r]^{P\otimes\id} & E\otimes_{\miq}\Omega_{\miq}^1 } $$
where $\iota$ is the natural inclusion.

We recall that $\nabla$ can be associated to a {\em connection form}
$A\in{\rm End}(W)\otimes$, so that $\nabla=\nabla_A=d+A$, see \cite[page 485]{FJ}.
Furthermore, there is a gauge in which $A=\Psi^{-1}d\Psi=-(d\Psi^{-1})\Psi$ for an
isomorphism $\Psi:W\otimes\miq\to E$ \cite[page 494]{FJ}.

\begin{proposition}
$\nabla$ is anti-self-dual.
\end{proposition}

The argument here is again very similar to the one in \cite[Proposition 14]{FJ};
we repeat it here for the sake of completeness.

\begin{proof}
Note that $F_{\nabla}=\nabla\nabla=PdPd$;
therefore we have:
\begin{eqnarray*}
F_{\nabla} & = & P \left( d ( \id_{\tW\otimes\miq} - 
\vec{\alpha}\Xi^{-1}\vec{\beta} ) d \right) =
P\left( d \vec{\alpha}\Xi^{-1} (d\vec{\beta}) \right) = \\
& = & P \left( (d\vec{\alpha})\Xi^{-1} (d\vec{\beta}) +
\vec{\alpha} d(\Xi^{-1} (d\vec{\beta}) ) \right) = \\
& = & P \left( (d\vec{\alpha})\Xi^{-1} (d\vec{\beta})  \right)
\end{eqnarray*}
for $P \vec{\alpha} d(\Xi^{-1} (d{\cal D}_{\rm I}) ) = 0$.
Since $\Xi^{-1}=(\beta_1\alpha_2)^{-1}\id_{\cpx^2}$, we conclude that 
$F_{\nabla}$
is proportional to $d\vec{\alpha}\wedge d\vec{\beta}$, as a 2-form.

It is then a straightforward calculation to show that each entry of
$d\vec{\alpha}\wedge d\vec{\beta}$ belongs to $\Omega^{2,-}_{\miq}$; indeed:
$$ d\vec{\alpha}\wedge d\vec{\beta} =
\left( \begin{array}{lr}
-dx_{11'} \ \ & \ \ -dx_{21'} \\
-dx_{12'} \ \ & \ \ -dx_{22'} \\ 0 & 0
\end{array} \right) \wedge
\left( \begin{array}{lrc}
dx_{22'} \ \ & \ \  -dx_{21'} & \ 0 \\
-dx_{12'} \ \ & \ \ dx_{11'} & \ 0
\end{array} \right) = $$
$$ = \left( \begin{array}{ccc}
-dx_{11'}dx_{22'}+dx_{21'}dx_{12'} \ \ &
\ \ dx_{11'}dx_{21'}-dx_{21'}dx_{11'} & \ 0 \\
-dx_{12'}dx_{22'}+dx_{22'}dx_{12'} \ \ &
\ \ dx_{12'}dx_{21'}-dx_{22'}dx_{11'} & \ 0 \\
0 & 0 & \ 0
\end{array} \right)$$
Applying the commutation relations (\ref{xnc2forms}), we obtain:
$$ d\vec{\alpha}\wedge d\vec{\beta} =
\left( \begin{array}{ccc}
-(dx_{11'}dx_{22'}+dx_{12'}dx_{21'}) &
2dx_{11'}dx_{21'} & \ 0 \\
-2dx_{12'}dx_{22'} &
dx_{11'}dx_{22'}+dx_{12'}dx_{21'} & \ 0 \\
0 & 0 & \ 0
\end{array} \right)$$
Comparison with (\ref{nc-asd2f}) completes the proof.
\end{proof}

\bigskip

The same procedure can be used to construct complex quantum instantons on 
$\mjq$; consider the maps:
$$ \xymatrix{ V\otimes\mjq \ar[r]^{\alpha_1}_{\alpha_2} &
\tW\otimes\mjq \ar[r]^{\beta_1}_{\beta_2} & V\otimes\mjq } $$
defined as follows:
$$ \begin{array}{lcr}
\alpha_1=\left( \begin{array}{c}
B_{11}\otimes\id - y_{22'} \\ B_{12}\otimes\id + \id\otimes y_{12'} \\ 
j_1\otimes\id
\end{array} \right) & {\rm and} &
\alpha_2=\left( \begin{array}{c}
B_{21}\otimes\id + y_{21'} \\ B_{22}\otimes\id - \id\otimes y_{11'} \\ 
j_2\otimes\id
\end{array} \right) \end{array} $$
$$ \beta_1=\left( \begin{array}{lcr}
-B_{12}\otimes\id - \id\otimes y_{12'} \ \ & \ \ B_{11}\otimes\id - 
\id\otimes y_{22'} \ \ & \ \ i_1\otimes\id
\end{array} \right) $$
$$ \beta_2=\left( \begin{array}{lcr}
-B_{22}\otimes\id + \id\otimes y_{11'} \ \ & \ \ B_{21}\otimes\id + 
\id\otimes y_{21'} \ \ & \ \ i_2\otimes\id
\end{array} \right) $$
Again, we set $\vec{\alpha}:(V\oplus V)\otimes\miq\to\tW\otimes\mjq$ and
$\vec{\beta}:\tW\otimes\miq\to (V\oplus V)\otimes\mjq$ given by:
$$ \vec{\alpha} = \left( \begin{array}{lr} \alpha_1 \ \ & \ \ \alpha_2 
\end{array} \right)
~~~ {\rm and} ~~~ \vec{\beta} = \left( \begin{array}{c} -\beta_2 \\ \beta_1 
\end{array} \right) ~ . $$
It follows that if $\vB$ is $\cpx$-regular, then $\Xi=\vec{\beta}\vec{\alpha}=\beta_1\alpha_2\id_{\cpx^2}$
is an isomorphism so that $E=\ker\vec{\beta}$ is a finitely generated free right $\mjq$-module;
an anti-self-dual connection is again produced via the projection formula.

The consistency map $\Gamma$ is obtained by restricting the obvious map
$\tW\otimes\miq[\det(x)^{-1}] \to \tW\otimes\mjq[\det(y)^{-1}]$, see (\ref{nc.glueing});
further details can be found in \cite{FJ}.

Finally, it is easy to see that $GL(V)$-equivalent complex ADHM data
will lead to gauge equivalent quantum instantons (see \cite{FJ}).
Our next result is the following:

\begin{theorem}
There is a well-defined map from the set of equivalence classes of $\cpx$-regular
solutions of the complex ADHM equations to the moduli space of gauge
equivalence classes of complex quantum instantons on $\mpq$.
\end{theorem}

By Proposition (\ref{proj}), $\cpx$-stability is a necessary condition for the
ADHM construction of instanton, possibly also sufficient. In that case, the domain
of the map in the theorem would be enlarged to the set of equivalence classes of
$\cpx$-stable solutions of the complex ADHM equations.


\subsection{Quantum Laplacian and admissibility} \label{qlap.sec}

In this final section we will relate the admissibility condition for sheaves on $\p^3$
and solutions of the Laplace equation in the quantum Minkowski space-time $\miq$,
thus extending the classical Penrose correspondence \cite{WW}. A $q$-deformation of the
Penrose transform with the quantum, rather than classical, twistor space has been
studied in \cite{SSSV} (see also the references in \cite{SSSV}).

We begin by constructing the quantum counterpart of the Laplace equation. Let
us define the quantum partial derivatives $\del_{rs'}$ by the relation:
\begin{equation} \label{del}
df = (\del_{11'}f)~dx_{11'} ~+~ (\del_{12'}f)~dx_{12'} ~+~
(\del_{21'}f)~dx_{21'} ~+~ (\del_{22'}f)~dx_{22'}
\end{equation}
where $f\in\miq$. Then the property $d^2=0$ implies the following commutation
relations:
\begin{equation} \label{delcommut1}
\del_{11'}\del_{21'} = \del_{21'}\del_{11'}, ~~~ \del_{12'}\del_{22'} = \del_{22'}\del_{12'},
\end{equation}
\begin{equation} \label{delcommut3}
[\del_{11'},\del_{22'}] + [\del_{12'},\del_{21'}] = 0
\end{equation}
\begin{equation} \label{delcommut2} \begin{array} {c}
\del_{11'}\del_{12'} = q^{-2}\del_{12'}\del_{11'},~~~
\del_{21'}\del_{22'} = q^{-2}\del_{22'}\del_{21'}, \\
\del_{12'}\del_{21'} = q^2 \del_{21'} \del_{12'}
\end{array} \end{equation}

We define the quantum Laplacian by:
\begin{equation} \label{q-lap}
\square_x = \del_{11'}\del_{22'} - \del_{21'}\del_{12'} = \del_{22'}\del_{11'} - \del_{12'}\del_{21'} ~ ,
\end{equation}
Alternatively, note that the quantum Laplacian can also be expressed via the Hodge star involution,
just as in the classical case:
\begin{equation} \label{q-lap.alt}
\square_x = * d * d ~ ,
\end{equation}
where the Hodge star on $0$- and $1$-forms is defined by the following (cf. \cite[Section 2.1]{FJ}):
\begin{equation} \label{H.0f}
*1 = q^{-1} dx_{11'}\wedge dx_{12'}\wedge dx_{21'}\wedge dx_{22'}
= q  dx_{22'}\wedge dx_{21'}\wedge dx_{12'}\wedge dx_{11'} 
\end{equation}
\begin{equation} \label{H.1f} \begin{array}{rcl}
\ast dx_{11'} & = & -\frac{1}{[2]} dx_{11'}\wedge dx_{12'}\wedge dx_{21'} \\
\ast dx_{12'} & = & -\frac{1}{[2]} dx_{12'}\wedge dx_{22'}\wedge dx_{11'} \\
\ast dx_{21'} & = & -\frac{1}{[2]} dx_{21'}\wedge dx_{11'}\wedge dx_{22'} \\
\ast dx_{22'} & = & -\frac{1}{[2]} dx_{22'}\wedge dx_{21'}\wedge dx_{12'}
\end{array} \end{equation}
where $[n]$ for $n\in\Z$, denotes the quantum integer $(q^n-q^{-n})/(q-q^{-1})$. 

We will now present a basis of solutions of the quantum Laplace equation
\begin{equation} \label{q-lap-eq}
\square_x f = 0 ~~ , ~~  f\in\miq
\end{equation}
in integral form, exactly as in the classical case \cite{WW}.

\begin{proposition} \label{q-lap-prop}
The following elements:
\begin{equation} \label{q-lap-soln}
X^l_{mn} = \frac{1}{2\pi i} \oint (x_{11'}s+x_{21'})^{l-m} (x_{12'}s+x_{22'})^{l+m} s^{n-l-1} ~ ds ~~,
\end{equation}
$$ -l \leq n,m \leq l ~ , ~ l\in\frac{1}{2}\Z_+ ~ , ~ n,m\equiv l ({\rm mod}~1) $$
form a basis of solutions of the quantum Laplace equation (\ref{q-lap-eq}),
where the integration variable $s$ commutes with the generators $x_{kl'}$.
\end{proposition}
\begin{proof}
First we will show that the above elements satisfy the quantum Laplace equation.
We note that since
$$ (x_{11'}s+x_{21'}) (x_{12'}s+x_{22'}) = (x_{12'}s+x_{22'}) (x_{11'}s+x_{21'}) $$
the expression (\ref{q-lap-soln}) does not depend on the order of the factors. Also, the
commutation relations (\ref{xnc1forms}) imply in particular:
\begin{equation} \label{d1}
(dx_{11'}s+dx_{21'}) (x_{11'}s+x_{21'}) = p^2 (x_{11'}s+x_{21'}) (dx_{11'}s+dx_{21'}) 
\end{equation}
\begin{equation} \label{d2}
(dx_{12'}s+dx_{22'}) (x_{12'}s+x_{22'}) = p^2 (x_{12'}s+x_{22'}) (dx_{12'}s+dx_{22'}) 
\end{equation}
for $p=q^{\pm1}$ and also
\begin{equation} \label{d3'}
(dx_{11'}s+dx_{21'}) (x_{12'}s+x_{22'}) = p^2 (x_{12'}s+x_{22'}) (dx_{11'}s+dx_{21'}) 
\end{equation}
for $p=q$ only and
\begin{equation} \label{d3''}
(dx_{12'}s+dx_{22'}) (x_{11'}s+x_{21'}) = p^2 (x_{11'}s+x_{21'}) (dx_{12'}s+dx_{22'}) 
\end{equation}
for $p=q^{-1}$ only. Therefore, in the computation of the differential $dX^l_{mn}$, the order
of factors has to be chosen accordingly. For $p=q$ we obtain:
$$ dX^l_{mn} = \frac{1}{2\pi i} \oint d\left( (x_{11'}s+x_{21'})^{l-m} (x_{12'}s+x_{22'})^{l+m} \right) s^{n-l-1} ~ ds = $$
$$ \left( \sum_{k=0}^{l+m-1} q^{2k} \right)  \oint (x_{11'}s+x_{21'})^{l-m} (x_{12'}s+x_{22'})^{l+m-1} (dx_{12'}s+dx_{22'}) s^{n-l-1} ~ ds + $$
$$ + \left( \sum_{k=l+m}^{2l-1} q^{2k} \right)  \oint (x_{11'}s+x_{21'})^{l-m-1} (x_{12'}s+x_{22'})^{l+m} (dx_{11'}s+dx_{21'}) s^{n-l-1} ~ ds $$
Similarly, for $p=q^{-1}$ we obtain:
$$ dX^l_{mn} = \frac{1}{2\pi i} \oint d\left(  (x_{12'}s+x_{22'})^{l+m} (x_{11'}s+x_{21'})^{l-m} \right) s^{n-l-1} ~ ds = $$
$$ = \left( \sum_{k=1}^{l-m-1} q^{-2k} \right)  \oint (x_{12'}s+x_{22'})^{l+m} (x_{11'}s+x_{21'})^{l-m-1} (dx_{11'}s+dx_{21'}) s^{n-l-1} ~ ds +$$
$$ + \left( \sum_{k=l-m}^{2l-1} q^{-2k} \right)  \oint (x_{12'}s+x_{22'})^{l+m-1} (x_{11'}s+x_{21'})^{l-m} (dx_{12'}s+dx_{22'}) s^{n-l-1} ~ ds $$
These yield the explicit expressions for the quantum partial derivatives, which can be
written uniformly for $p=q^{\pm1}$ as follows:
\begin{equation}\label{partials}\begin{array}{c}
\del_{11'} X^l_{mn} = p^{2l-1}q^{m+l} ~ [l-m] ~ X^{l-\hf}_{(n+\hf)(m+\hf)} \\
\del_{12'} X^l_{mn} = p^{2l-1}q^{m-l} ~ [l+m] ~ X^{l-\hf}_{(n+\hf)(m-\hf)} \\
\del_{21'} X^l_{mn} = p^{2l-1}q^{m+l} ~ [l-m] ~ X^{l-\hf}_{(n-\hf)(m+\hf)} \\
\del_{22'} X^l_{mn} = p^{2l-1}q^{m-l} ~ [l+m] ~ X^{l-\hf}_{(n-\hf)(m-\hf)}
\end{array}\end{equation}
The relations (\ref{partials}) immediately imply:
\begin{equation}\label{q-soln} \square_x X^l_{nm} =0 \end{equation}

Finally, we need to argue that any solution of the quantum Laplace equations (\ref{q-lap-eq})
is a (complex) linear combination of the ones above. In fact, as in the classical case, the
elements:
\begin{equation}\label{detX}
(\det(x))^k X^l_{nm} ~ , ~ k\in\Z_+  ~ , ~ l\in\hf\Z_+ ~ , ~ -l \leq n,m\leq l 
\end{equation}
$$ {\rm where} ~~ \det(x)= x_{11'}x_{22'}-x_{12'}x_{22'} = x_{22'}x_{11'} - x_{21'}x_{12'} $$
form a basis of $\miq$, since they are linearly independent (even for $q=1$),
and the number of elements of fixed degree in (\ref{detX}) is the same as the
number of ordered monomials  (\ref{basis}) on the four variables $x_{rs'}$,
which also compose a basis of $\miq$.

On the other hand, for $q=1$, the elements (\ref{detX}) with $k=0$ form a basis of
solutions of the Laplace equations, known as harmonic polynomials, and this space
cannot increase for generic or formal parameter $q$.
\end{proof} 

In fact, one can prove an explicit quantum analog of the spectral decomposition of the Laplace
operator, which contains, as a special case, the statement of the Proposition 34.

\begin{proposition} \label{laptil}
The basis (\ref{detX}) consists of eigenfunctions of the operator
\begin{equation} \label{q-laptil}
\tilde{\square}_x = \det(x) \cdot \square_x
\end{equation}
with eigenvalues $p^{2k+2l-3}[k][k+2l+1]$.
\end{proposition}

The proof requires an elementary quantum calculus, which will be given in Appendix \ref{app-b}.

Now we would like to reinterpret Proposition \ref{q-lap-prop} in terms of the sheaf cohomology of
$\p^{\rm I}=\p^3\setminus\ell_\infty$. We will consider a covering of $\p^{\rm I}$ by its two simply
connected patches:
$$ \p^{\rm I}_{(1)}=\{ [x:y:z:w]\in\p^3 ~ | ~ z\neq0 \} $$
$$ \p^{\rm I}_{(2)}=\{  [x:y:z:w]\in\p^3 ~ | ~ w\neq0 \} $$
Then the elements of $H^1(\p^{\rm I},{\cal O}_{\p^{\rm I}}(-2))$ can be represented by transition
functions which are rational in $x,y,z,w$ of degree $-2$, with singularities along the hyperplanes
$\{z=0\}$ and $\{w=0\}$. A natural basis of such functions is given by the Laurent monomials
\begin{equation}\label{cech}
\frac{x^{l-m}~y^{l+m}}{z^{l-n+1}~w^{l+n+1}} ~ , ~ l\in\hf\Z_+ ~ , ~ -l\leq n,m \leq l ~.
\end{equation}

The quantum Penrose transform assigns to an element of this basis the quantum harmonic
polynomials $X^l_{nm}$ via the formula (\ref{q-lap-soln}). Thus we can restate Proposition
\ref{q-lap-prop} as follows.

\begin{proposition} \label{PT1}
There is an isomorphism
$$ H^1(\p^{\rm I},{\cal O}_{\p^{\rm I}}(-2)) \stackrel{\simeq}{\longrightarrow} \ker\square_x $$
given by the integral formula
$$ f(x,y,z,w) \mapsto  \frac{1}{2\pi i} \oint f(x_{11'}s+x_{21'} , x_{12'}s+x_{22'} , s , 1) ~ ds $$
where $f$ represents a cocycle in $H^1(\p^{\rm I},{\cal O}_{\p^{\rm I}}(-2))$, that is, a linear
combination of the terms in (\ref{cech}).
\end{proposition}

Next we will extend this isomorphisms to more general vector bundles over $\p^{\rm I}$. More
precisely, if a vector bundle $\cE^{\rm I} \to \p^{\rm I}$ is admissible (i.e. it is the restriction of
an admissible vector bundle $\cE\to\p^3$ to $\p^{\rm I}$), then a class in $H^1(\p^{\rm I},\cE^{\rm I}(-2))$
can be represented by a cocycle $\vec{f}=(f_1,\cdots,f_r)$, where $r$ is the rank of  $\cE^{\rm I}$,
and each $f_k$ is a linear combination of the terms in (\ref{cech}).

On the other hand, given a complex quantum instanton $(E,\nabla_A)$ on $\miq$ we define the quantum
coupled Laplacian by generalizing (\ref{q-lap.alt}):
\begin{equation} \label{q-c-lap}
\square^A_x = * \nabla_A * \nabla_A ~~.
\end{equation}

To any given admissible vector bundle $\cE^{\rm I} \to \p^{\rm I}$ we can associate a regular ADHM
data $\vB$, which can then be used to construct a complex quantum instanton $(E,\nabla_A)$ on $\miq$.
For $\cE^{\rm I}$ and $(E,\nabla_A)$ related as above, our next theorem generalizes the classical Penrose
correspondence between $H^1(\p^{\rm I},\cE^{\rm I}(-2))$ and the set of solutions of the coupled quantum
Laplace equation. 

\begin{theorem} \label{PT2}
There is an isomorphism
$$ H^1(\p^{\rm I},\cE^{\rm I}(-2)) \stackrel{\simeq}{\longrightarrow} \ker\square_x^A $$
given by the integral formula
$$ f(x,y,z,w) \mapsto  \varphi = \frac{1}{2\pi i} \oint \Psi^{-1} \vec{f}(x_{11'}s+x_{21'} , x_{12'}s+x_{22'} , s , 1) ~ ds $$
\end{theorem}
\begin{proof}
First recall that there is a gauge such that $\nabla_A = d - (d\Psi^{-1})\Psi$.
We then have:
\begin{eqnarray*}
\nabla_A \varphi & = & \frac{1}{2\pi i} \oint (d - (d\Psi^{-1})\Psi) \Psi^{-1}
\vec{f}(x_{11'}s+x_{21'} , x_{12'}s+x_{22'} , s , 1) ~ ds = \\
& = &  \frac{1}{2\pi i} \oint \left( (d\Psi^{-1})\vec{f} + \Psi^{-1}d\vec{f} - (d\Psi^{-1})\vec{f}  \right) ~ ds =
\frac{1}{2\pi i} \oint \left( \Psi^{-1}d\vec{f} \right)~ ds 
\end{eqnarray*}
By the same token, it follows that: 
\begin{eqnarray*}
\square_x^A \varphi & = & * \frac{1}{2\pi i} \oint (d - (d\Psi^{-1})\Psi) \Psi^{-1}(*d\vec{f}) ~ ds = \\ 
& = & * \frac{1}{2\pi i} \oint \Psi^{-1}(d*d\vec{f}) ~ ds = \frac{1}{2\pi i} \oint \Psi^{-1}(\square_x\vec{f}) ~ ds = 0 
\end{eqnarray*}
since $\square_x\vec{f}=0$ by Proposition \ref{PT1}.
\end{proof}
  
We now consider global solutions of the quantum Laplace equation on $\mpq$.
To do that, we have to check the consistency of the solutions in $\miq$ and $\mjq$.

We introduce the following elements of $\mjq$:
\begin{equation} \label{soln-y}
Y^l_{nm} = \frac{1}{2\pi i} \oint (y_{11'}t+y_{12'})^{l-n} (y_{21'}t+y_{22'})^{l+n} t^{m-l-1} ~ dt ~~,
\end{equation}
$$ -l \leq n,m \leq l ~ , ~ l\in\frac{1}{2}\Z_+ ~ , ~ n,m\equiv l ({\rm mod}~1) ~ .$$
Then one defines the quantum partial derivatives and the quantum Laplacian
$\square_y$ in $\mjq$ as above. The counterparts of Propositions \ref{q-lap-prop} and \ref{laptil}
hold with the eigenvalues of $\tilde{\square}_y$ being equal to $p^{-2k-2l+3}[k][k+2l+1]$ on the
basis elements $(\det(y))^kY^l_{nm}$, where $\det(y)$ was defined in Section \ref{qinst1}

Including the negative powers of $\det(y)$ in the above basis and the negative powers of
$\det(x)$ in the basis (\ref{detX}), we obtain two bases of $\mijq$. We can also extend the
quantum Laplacians $\square_x$ and $\square_y$ in $\mijq$ so that Proposition \ref{laptil} and
its counterpart in the $y$ generators hold. The comparison between the two systems of coordinates
gives essentially the same result as in the classical $q=1$ case.

\begin{proposition} \label{oast}
In $\mijq$, one has:
\begin{itemize}
\item  $(\det(x))^kX^l_{nm}$ is proportional to $(\det(y))^{-k-2l}Y^l_{nm}$;
\item $\tilde{\square}_x = p^{-8} (\det(y)) \tilde{\square}_y (\det(y))^{-1}$.
\end{itemize} \end{proposition}

\begin{proof}
The first statement follows from the definitions (\ref{q-lap-soln}) and (\ref{soln-y}); for the
details, see Appendix \ref{app-b}.

The second identity follows from the comparison of the eigenvalues of the two operators
in the proportional bases $(\det(x))^kX^l_{nm}$ and $(\det(y))^{-k-2l}Y^l_{nm}$, as
computed in Proposition \ref{laptil}. It can also be verified directly by relating the partial
derivatives with respect to $x_{kl'}$ and $y_{kl'}$.
\end{proof}

Thus, as in the classical case, we conclude that there are no consistent solutions to the
scalar quantum Laplace equation in the (compactified) quantum Minkowski space-time
$\mpq$. With Proposition \ref{PT1} in mind, this non-existence statement corresponds
to the fact that $H^1(\p^3,{\cal O}(-2))=0$.

More generally, let $\cE$ be an admissible vector bundle over $\p^3$, to which we can
associate, via the intermediate $\cpx$-regular ADHM datum, a complex quantum instanton
$(E_{\rm I},\nabla_{\rm I};E_{\rm J},\nabla_{\rm J})$ over $\mpq$. Taking the non-existence
of consistent solutions to the scalar quantum Laplace equation to the context of Theorem
\ref{PT2}, we conclude that there are no consistent solutions to the coupled quantum
Laplace equation.This in turn corresponds to the vanishing of $H^1(\p^3,\cE(-2))$.

In a future paper, we plan to establish the reverse correspondence: given a complex
quantum instanton over $\mpq$, we will associate directly an admissible vector bundle
over $\p^3$. This will generalize the celebrated Penrose-Ward correspondence.


\appendix
\section{Moduli space of stable ADHM data}\label{ap-nak}

Here we recall the proof that ${\cal M}(r,c)$ is a smooth complex 
manifold of dimension $2rc$; some of the arguments were relevant
in Section \ref{adhm}, for the proof of Theorem \ref{main1}. Our
arguments are inspired by \cite{H,N2,V}.

Let $V$ and $W$ be complex vector spaces, with dimensions $c$
and $r$, respectively, and set $\tW=V\oplus V\oplus W$. Define also:
$$ \mathbf{B} = {\rm Hom}(V,V)\oplus{\rm Hom}(V,V)\oplus
{\rm Hom}(W,V)\oplus{\rm Hom}(V,W) $$
A point $(B_{k},i,j)\in\mathbf{B}$ ($k=1,2$) is called a {\it ADHM datum}.
As mentioned above, the groups $GL(V)$ and $GL(W)$ act on
$\mathbf{B}$ in the following way:
\begin{equation} \label{V-action}
g\cdot (B_k,i,j) = (gB_kg^{-1},gi,jg^{-1}) ~ , ~~~ g\in GL(V)
\end{equation}
\begin{equation} \label{W-action}
g\cdot (B_k,i,j) = (B_k,ig^{-1},gj) ~ , ~~~ g\in GL(W)
\end{equation}

\begin{theorem} \label{nakvar}
${\cal M}(r,c)$ is a smooth, complex manifold of dimension $2rc$.
\end{theorem}

Indeed, it is known that ${\cal M}(r,c)$ is non-empty for all $r,c\geq1$. 
Furthermore, it can be shown that ${\cal M}(r,c)$ is a simply-connected
quasi-projective algebraic variety \cite{V} and that it admits a complete
hyperk\"ahler metric \cite{N2}. The strategy of the proof goes as follows;
considering the map:
\begin{eqnarray*}
& \mu : \mathbf{B}^{\rm st}\to{\rm Hom}(V,V) & \\ & \mu(B_1,B_2,i,j)=[B_1,B_2] +ij &
\end{eqnarray*}
where $\mathbf{B}^{\rm st}$ is the open subset of stable ADHM data. We
first show that $GL(V)$ acts freely in $\mathbf{B}^{\rm st}$, and that the
action has a closed graph; we then show that $\mu^{-1}(0)$ is indeed a
complex manifold of dimension $2rc+c^2$; the desired result follows from
general theory (see for instance the {\em closed graph lemma} in \cite{H}
and the references therein).

\begin{proposition} \label{stability}
$(B_1,B_2,i,j)$ is stable if and only if:
\begin{enumerate}
\item $(B_1,B_2,i,j)$ is not fixed of the $GL(V)$ action;
\item if $X\in{\rm Hom}(V,V)$ satisfies
$ [ B_1, X ] =   [ B_2, X ] =  X i = 0$, then $X=0$.
\end{enumerate}\end{proposition}
\begin{proof}
Suppose that $(B_1,B_2,i,j)$ is fixed by some $g\neq\id_V\in GL(V)$, so 
that, $gB_kg^{-1}=B_k$ ($k=1,2$) and $gi=i$. The former implies that
$\ker(g-\id_V)$ is $B_k$ invariant, while the latter implies that
$\ker(g-\id_V)\subset{\rm Im}i$, thus contradicting stability.

For the second statement, $Xi=0$ implies that $i(W) \subset \ker X$, while 
$[ B_k, X ]  =  0$ implies that $\ker X$ is $B_k$-invariant. Stability then
implies that $X=0$.
\end{proof}

\begin{proposition} \label{proper}
The action (\ref{V-action}) has a closed graph, i.e. the set
$$ \Gamma = \{ (X,Y) \in \mathbf{B}^{\rm st}\times\mathbf{B}^{\rm st} ~ | ~
Y=g \cdot X ~ {\rm for~some} ~ g\in GL(V) \} $$
is closed in $\mathbf{B}^{\rm st}\times\mathbf{B}^{\rm st}$. In other words,
$GL(V)$ acts properly in $\mathbf{B}^{\rm st}$.
\end{proposition}
\begin{proof}
Let $\{X_k\}$ be a sequence in
$\mathbf{B}^{\rm st}$, while $\{g_k\}$ denotes a sequence in $GL(V)$;
assuming that:
$$ \lim_{k\to\infty} X_k=X_\infty ~~~ {\rm and} ~~~
\lim_{k\to\infty} g_k\cdot X_k=Y_\infty ~,$$
we must show that $Y_\infty =g_\infty\cdot X_\infty$ for some $g_\infty\in GL(V)$,
or equivalently that the sequence $\{g_k\}$ converges to $g_\infty\in GL(V)$.

Indeed, for any given $X=(B_1,B_2,i,j)\in\mathbf{B}^{\rm st}$ we consider
the map $R (X) : W^{\oplus c^2} \longrightarrow V$ given by ($1\leq m,n \leq c-1$):
$$ R(X) = i \oplus \cdots \oplus B_1^mB_2^n i \oplus \cdots \oplus B_1^{c-1}B_2^{c-1} i ~ .$$ 
Note that $gR(X)=R(g\cdot X)$ for any $g\in GL(V)$.

Furthermore, $R(X)$ is surjective if and only if $X$ is stable. Indeed, if
$X=(B_1,B_2,i,j)$ is not stable, then there is $v\in V^*$ such that
$B_1^*v=\lambda_1v$, $B_2^*v=\lambda_2v$ and $i^*v=0$; hence
$R^*v=0$ so that $R$ is not surjective. Conversely, if $R$ is not surjective,
then $S={\rm Im}~R$ is a proper subspace of $V$; clearly, $S$ is $B_1$
and $B_2$ invariant, and $i(W)\subset S$, hence $X$ is not stable. 

The sequence of maps $R(X_k)$ converges to $R(X_\infty)$; thus, there
is a sequence of maps $T_k\in{\rm Hom}(V,W^{\oplus c^2})$ converging
to a map $T_\infty\in{\rm Hom}(V,W^{\oplus c^2})$ such that:
$$ W^{\oplus c^2} = \ker R(X_k) \oplus {\rm Im}~T_k = \ker R(X_\infty) \oplus {\rm Im}~T_\infty $$ 
It then follows that $R(X_k)T_k$ and $R(X_\infty)T_\infty$ are invertible
as operators on $V$.

Now set $g_\infty = R(Y_\infty)T_\infty [R(X_\infty)T_\infty]^{-1} \in GL(V)$. Thus:
$$ g_k = g_k [R(X_k)T_k] [R(X_k)T_k]^{-1} = [R(g_k X_k)T_k] [R(X_k)T_k]^{-1} $$
and $g_k$ converges to $g_\infty$.
\end{proof}

\begin{proposition} \label{DBsurj}
$\vB=(B_1,B_2,i,j)$ is stable if and only if the derivative map
$D_{\vB}\mu:\mathbf{B}\to{\rm Hom}(V,V)$ is surjective.
\end{proposition}
This means that $0$ is a regular value of the map $\mu$, hence $\mu^{-1}(0)$
is a smooth complex manifold of dimension $2rc+c^2$.
\begin{proof}
Taking $(b_1,b_2,c,d)\in\mathbf{B}$, the derivative map is given by:
$$ D_{\vB}\mu(b_1,b_2,c,d) = [b_1,B_2] + [B_1,b_2] + id + cj $$
Let $X\in{\rm Hom}(V,V)$ be orthogonal to the image of $D_{\vB}\mu$, that is
$$ {\rm Tr}(D_B\mu(b_1,b_2,c,d)X^\dagger)=0 ~, ~~ \forall (b_1,b_2,c,d) ~ .$$
Then in particular:
$$ {\rm Tr}([b_1,B_2]X^\dagger)={\rm Tr}(b_1[X^\dagger,B_2])=0 ~~
\forall b_1 $$
$$ {\rm Tr}([B_1,b_2]X^\dagger)={\rm Tr}([X^\dagger,B_1]b_2)=0 ~~
\forall b_2 $$
$$ {\rm Tr}(idX^\dagger) = {\rm Tr}(X^\dagger id)=0 ~~ \forall d $$
Hence $[X^\dagger,B_1]=[X^\dagger,B_1]=X^\dagger i=0$, so $X=0$ by 
Proposition \ref{stability}.
\end{proof}

Since $GL(V)$ acts freely and properly on the smooth manifold $\mu^{-1}(0)$,
this completes the proof of Theorem \ref{nakvar}. To conclude this section,
we also remark upon the following statements, in which by {\em irregular} we
mean neither stable nor costable:

\begin{proposition} \label{r/ir}
Every solution of (\ref{adhm1}) and (\ref{adhm2}) is:
\begin{enumerate}
\item stable, if $\xi>0$;
\item costable, if $\xi<0$;
\item either regular or irregular, if $\xi=0$.
\end{enumerate}\end{proposition}
\begin{proof}
For the first statement, if $(B_1,B_2,i,j)$ is not stable, then by
duality on $V$ there is a proper subspace $S^\perp\subset V$ such that
$B_k^\dagger(S^\perp)\subset S^\perp$ and $S^\perp\subset\ker i^\dagger$.
So restricting (\ref{adhm2}) to $S^\perp$ and taking the trace,
we conclude that
$$ {\rm Tr}(ii^\dagger|_{S^\perp}) = \xi\cdot\dim S + {\rm Tr}(j^\dagger 
j|_{S^\perp}) > 0 $$
which yields a contradiction.
The proof of the second statement is similar, while the third
statement can be found at \cite[Lemma 2]{FJ}.
\end{proof}

It is interesting to compare the third part of Proposition \ref{r/ir} with
Remark \ref{st.not.reg}: the complex equations are a much more flexible
than the real ones. 

\begin{proposition} \label{r=1case} \cite[p. 24]{N2}.
Let $r=1$. Every stable solution of (\ref{adhm1}) has $j=0$. In particular,
there are no regular solutions for $r=1$ and $\xi=0$.
\end{proposition}


\section{Cohomological calculations} \label{ap1}

We collect here the proofs for various facts used in Section \ref{p3}.

\begin{proposition} \label{ap-prop1}
Let $\cE$ be an admissible torsion-free sheaf over $\p^3$ with
${\rm ch}(\cE) = r - c [H]^2$ and such that
$\cE|_{\ell_\infty}={\cal O}_{\ell_\infty}^{\oplus r}$. The 
following hold:
\begin{enumerate}
\item $h^1(\p^3,\cE(-1))=-\chi(\cE(-1))=c$;
\item $h^1(\p^3,\cE\otimes\Omega^1_{\p^3})=-\chi(\cE\otimes\Omega^1_{\p^3})=c+2r$;
\item $h^1(\p^3,\cE\otimes\Omega^2_{\p^3}(1))=-\chi(\cE\otimes\Omega^2_{\p^3}(1))=c$;
\item $H^1(\p^3,\cE\otimes\Omega^2_{\p^3}(1))\simeq H^1(\p^3,\cE(-1))$.
\end{enumerate}\end{proposition}
\begin{proof}
Let us first spell out the admissibility condition more precisely:
\begin{equation} \begin{array}{ccc}
H^0(\p^3,\cE(k))=0,~~\forall k\leq-1 & \ \ & H^1(\p^3,\cE(k))=0,~~\forall 
k\leq-2 \\
H^2(\p^3,\cE(k))=0,~~\forall k\leq-2 & \ \ & H^3(\p^3,\cE(k))=0,~~\forall 
k\leq-3
\end{array}\end{equation}
The first statement then follows immediately from admissibility, and it only
remains for us to show that $\chi(\cE(-1))=-c$. Indeed, note that:
$$ {\rm ch}(\cE(-1)) = r - r h + \left( \frac{r}{2}+c\right) h^2 +
\left( -\frac{r}{6}+c\right) h^3 $$
$$ {\rm td}(\p^3) = 1 + 2h + \frac{22}{12}h^2 + h^3 $$
Hence it follows:
$$ \chi(\cE(-1)) = \int_{\p^3} {\rm ch}(\cE(-1)){\rm td}(\p^3) = -c $$

Now consider the Euler sequence for 1-forms:
\begin{equation} \label{e1f}
0\to \Omega_{\p^3}^1 \to \bigoplus_4 \op3(-1) \to \op3 \to 0
\end{equation}
from which we conclude that:
$$ {\rm ch}(\Omega_{\p^3}^1) = 3 - 4h + 2h^2 + \frac{2}{3}h^3 $$
$$ {\rm ch}(\cE\otimes\Omega_{\p^3}^1) = 3r + 4r h - (3c-2r)h^2 -
\left( \frac{2r}{3}+4c\right) h^3 $$
Using Riemann-Roch again we obtain:
$$ \chi(\cE\otimes\Omega_{\p^3}^1) =
\int_{\p^3} {\rm ch}(\cE\otimes\Omega_{\p^3}^1){\rm td}(\p^3) = -c-2r $$
Tensoring (\ref{e1f}) by $\cE$ we obtain the exact sequence:
$$ {\rm Tor}^1(\cE,\op3)\to \cE\otimes\Omega_{\p^3}^1 \to \bigoplus_4 
\cE(-1) \to \cE \to 0 $$
But the first term vanishes because $\op3$ is a locally-free sheaf. 
Therefore we have:
\begin{equation} \label{ee1f}
0\to \cE\otimes\Omega_{\p^3}^1 \to \bigoplus_4 \cE(-1) \to \cE \to 0
\end{equation}
At the level of cohomology, one obtains:
$$ 0 \to H^0(\p^3,\cE\otimes\Omega^1_{\p^3}) \to \bigoplus_4 
H^0(\p^3,\cE(-1)) $$
and since $H^0(\p^3,\cE(-1))=0$, it follows that 
$H^0(\p^3,\cE\otimes\Omega^1_{\p^3})=0$.
Moreover, (\ref{ee1f}) also implies that:
$$ H^2(\p^3,\cE) \to H^3(\p^3,\cE\otimes\Omega^1_{\p^3}) \to \bigoplus_4 
H^3(\p^3,\cE(-1)) $$
Since the first and third groups vanish by admissibility, we obtain
$H^3(\p^3,\cE\otimes\Omega^1_{\p^3})=0$.

The Euler sequence for 3-forms is given by:
\begin{equation} \label{e3f}
0\to \Omega_{\p^3}^3 \to \bigoplus_4 \op3(-3) \to \Omega_{\p^3}^2 \to 0
\end{equation}
Recalling that $\Omega_{\p^3}^3=\op3(-4)$ and tensoring (\ref{e3f}) by 
$\cE(1)$,
we obtain:
\begin{equation} \label{ee3f}
0\to \cE(-3) \to \bigoplus_4 \cE(-2) \to \cE\otimes\Omega_{\p^3}^2(1) \to 0
\end{equation}
Since $H^2(\p^3,\cE(-2))=0$ for all $p$, it follows from the cohomology 
sequence
associated with (\ref{ee3f}) that:
\begin{equation} \label{iso1}
H^p(\p^3,\cE\otimes\Omega_{\p^3}^2(1)) \simeq H^{p+1}(\p^3,\cE(-3))
\end{equation}
thus $H^p(\p^3,\cE\otimes\Omega_{\p^3}^2(1))=0$ for $p=0,2,3$ by 
admissibility.
The sequence (\ref{ee3f}) can also be used to compute the Chern character
of $\cE\otimes\Omega_{\p^3}^2(1)$; indeed,
\begin{eqnarray*}
{\rm ch}(\cE\otimes\Omega_{\p^3}^2(1)) & = & 4{\rm ch}(\cE(-2)) - {\rm 
ch}(\cE(-3)) = \\
& = & 3r - 5rh + 
\left(\frac{7r}{2}-3c\right)h^2+\left(5c-\frac{5r}{6}\right)h^3
\end{eqnarray*}
It then follows that:
$$ \chi(\cE\otimes\Omega_{\p^3}^2(1)) =
\int_{\p^3} {\rm ch}(\cE\otimes\Omega_{\p^3}^2(1)){\rm td}(\p^3) = - c $$
what completes the proof of the third statement.

To complete the proof of the second statement, it only remains for us to show
that $H^2(\p^3,\cE\otimes\Omega^1_{\p^3})=0$. Tensoring the Euler sequence
for 2-forms
\begin{equation} \label{e2f}
0\to \Omega_{\p^3}^2 \to \bigoplus_6 \op3(-2) \to \Omega_{\p^3}^1 \to 0
\end{equation}
by $\cE$ we obtain:
\begin{equation} \label{ee2f}
0\to \cE\otimes\Omega_{\p^3}^2 \to \bigoplus_6 \cE(-2) \to 
\cE\otimes\Omega_{\p^3}^1 \to 0
\end{equation}
The associated cohomology sequence yields:
$$ \bigoplus_6 H^{2}(\p^3,\cE(-2)) \to H^2(\p^3,\cE\otimes\Omega_{\p^3}^1) 
\to H^3(\p^3,\cE\otimes\Omega_{\p^3}^2) \to \bigoplus_6 H^{3}(\p^3,\cE(-2)) $$
Admissibility implies that the first and last groups vanish, hence
$H^2(\p^3,\cE\otimes\Omega_{\p^3}^1)\simeq 
H^3(\p^3,\cE\otimes\Omega_{\p^3}^2)$.
Now tensoring (\ref{e3f}) by $\cE$ we get:
\begin{equation} \label{ee3f'}
0\to \cE(-4) \to \bigoplus_4 \cE(-3) \to \cE\otimes\Omega_{\p^3}^2 \to 0
\end{equation}
we conclude that $H^3(\p^3,\cE\otimes\Omega_{\p^3}^2)=0$ since 
$H^3(\p^3,\cE(-3))=0$ by admissibility.

Finally, let $\wp$ be a plane containing $\ell_\infty$, so that the restriction
$\cE|_{\wp}$ yields a torsion-free sheaf on $\wp$ which is trivial at 
$\ell_\infty$. Consider the sequence:
$$ 0 \to \cE(-p-1) \to \cE(-p) \to \cE(-p)|_{\wp} \to 0 $$
Setting $p=-2$, we conclude that $H^1(\p^3,\cE(-2)|_{\wp}) \simeq 
H^2(\p^3,\cE(-3))$,
since $H^1(\p^3,\cE(-2))=H^2(\p^3,\cE(-2))=0$ by admissibility. Then setting 
$p=-1$,
we get that $H^1(\p^3,\cE(-1)|_{\wp}) \simeq H^1(\p^3,\cE(-1))$ for the same 
reason.
Together with (\ref{iso1}), we have obtained the identifications
$$ H^1(\p^3,\cE\otimes\Omega_{\p^3}^2(1)) \simeq H^{2}(\p^3,\cE(-3)) \simeq $$
$$ \simeq H^1(\p^3,\cE(-2)|_{\wp}) \simeq H^1(\p^3,\cE(-1)|_{\wp})
\simeq H^1(\p^3,\cE(-1)) $$
where the third identification follows from \cite[page 20]{N2}. This 
completes the proof of the fourth statement.
\end{proof}

\begin{proposition} \label{ap-prop2}
Consider the following exact sequence of sheaves on a
regular algebraic variety $V$ of dimension 3:
\begin{equation} \label{sqc.x}
0 \to {\cal A} \stackrel{\mu}{\longrightarrow} {\cal B}
\longrightarrow {\cal C} \to 0
\end{equation}
where $\cal A$ and $\cal B$ are locally-free. Then:
\begin{enumerate}
\item $\cal C$ is torsion-free if and only if the localized map
$\mu_X:{\cal A}_X\to{\cal B}_X$ is injective away from a subset
of codimension 2;
\item $\cal C$ is reflexive if and only if the localized map
$\mu_X:{\cal A}_X\to{\cal B}_X$ is injective away from a subset
of codimension 3.
\item $\cal C$ is locally-free if and only if the localized map
$\mu_X:{\cal A}_X\to{\cal B}_X$ is injective for all $X\in V$.
\end{enumerate} \end{proposition}

\begin{proof}
Dualizing the sequence (\ref{sqc.x}) we obtain:
\begin{equation} \label{sqc dual}
0 \to {\cal C}^* \longrightarrow {\cal B}^* 
\stackrel{\mu^*}{\longrightarrow} {\cal A}^*
\longrightarrow {\rm Ext}^1({\cal C},\op3) \to 0
\end{equation}
It follows that ${\rm Ext}^p({\cal C},\op3)=0$ for $p=2,3$, and
$$ I = {\rm supp}\left({\rm Ext}^1({\cal C},\op3)\right) =
\{ X\in V ~|~ \mu_X~{\rm is~not~injective}\} $$
So it is now enough to argue that $\cal C$ is torsion-free if and only if
$\dim I=1$ and that $\cal C$ is reflexive if and only if $\dim I=0$.
The third statement is clear.

Recall that the $m^{\rm th}$-singularity set of a coherent sheaf $\cal F$
is given by:
$$ S_m({\cal F}) = \{ X\in\p^3 ~|~ dh({\cal F}_x) \geq 3-m \} $$
where $dh({\cal F}_x)$ stands for the homological dimension of
${\cal F}_x$ as an ${\cal O}_x$-module:
$$ dh({\cal F}_x) = d ~~~ \Longleftrightarrow ~~~
\left\{ \begin{array}{l}
{\rm Ext}^d_{{\cal O}_x}({\cal F}_x,{\cal O}_x) \neq 0 \\
{\rm Ext}^p_{{\cal O}_x}({\cal F}_x,{\cal O}_x) = 0 ~~ \forall p>d
\end{array} \right. $$

In the case at hand, we have that $dh({\cal F}_x) = 1$ if $X\in I$,
and $dh({\cal F}_x) = 0$ if $X\notin I$. Therefore
$S_0({\cal C})=S_1({\cal C})=\emptyset$, while $S_2({\cal C})=I$.
It follows that \cite[Proposition 1.20]{ST} :
\begin{itemize}
\item if $\dim I = 1$, then $\dim S_m({\cal C})\leq m-1$ for all $m<3$,
hence $\cal C$ is a locally 1$^{\rm st}$-syzygy sheaf;
\item  if $\dim I = 0$, then $\dim S_m({\cal C})\leq m-2$ for all $m<3$,
hence $\cal C$ is a locally 2$^{\rm nd}$-syzygy sheaf.
\end{itemize}
The desired statements follow from the observation that $\cal C$ is
torsion-free if and only if it is a locally 1$^{\rm st}$-syzygy sheaf, 
while $\cal C$ is reflexive if and only if it is a locally
2$^{\rm nd}$-syzygy sheaf \cite[page 148-149]{OSS}. 
\end{proof}


\section{Quantum space-time calculations} \label{app-b}

We collect here the proofs of various facts used in Section \ref{qlap.sec}.

First, we will derive a few formulas for the differential forms on quantum
Minkowski space-time $\miq$. The commutation relations between $x_{ij'}$
and $dx_{kl'}$ imply:
\begin{equation} \label{dx-det} \begin{array}{rcl}
dx_{11'} ~ \det(x) & = & p^2  \det(x) ~ dx_{11'} \\
dx_{12'} ~ \det(x) & = & p^2q^{-2}  \det(x) ~ dx_{12'} \\
dx_{21'} ~ \det(x) & = & p^2q^2  \det(x) ~ dx_{21'} \\
dx_{22'} ~ \det(x) & = & p^2  \det(x) ~ dx_{22'}
\end{array} \end{equation}
We also obtain:
\begin{eqnarray} 
\label{ddet} d(\det(x)) & = & dx_{11'}~x_{22'} + x_{11'}~dx_{22'} - dx_{12'}~x_{21'} - x_{12'}~dx_{21'} = \\
\nonumber & = & p^{-1}q \left( x_{11'}dx_{22'} - x_{12'} d x_{21'} \right) +
p^{-1}q^{-1} \left(x_{22'}dx_{11'} - x_{21'}dx_{12'}\right)
\end{eqnarray}
Applying Leibnitz rule we obtain:
$$ d (f\det(x)) = \left( \sum_{kl'} (\partial_{kl'}f) dx_{kl'}\right)\det(x) + fd(\det(x)) $$
combining this with (\ref{dx-det}) and  (\ref{ddet}) we have: 
\begin{equation} \label{del-det} \begin{array}{rcl}
\partial_{11'} (f\det(x)) & = & p^2  (\partial_{11'}f) \det(x) +  p^{-1}q^{-1} fx_{22'} \\
\partial_{12'} (f\det(x)) & = & p^2q^{-2} (\partial_{12'}f) \det(x) -  p^{-1}q^{-1} fx_{21'} \\
\partial_{21'} (f\det(x)) & = & p^2q^2  (\partial_{21'}f) \det(x) -  p^{-1}q fx_{12'} \\
\partial_{22'} (f\det(x)) & = & p^2 (\partial_{22'}f) \det(x) +  p^{-1}q fx_{11'} 
\end{array} \end{equation}

We introduce an operator:
\begin{equation} \label{delta-op}
\Delta_x f = (\partial_{11'}f) x_{11'} +  (\partial_{12'}f) x_{12'}  +
(\partial_{21'}f) x_{21'} + (\partial_{22'}f) x_{22'}
\end{equation}
and let $D_x$ denote the operator of multiplication by $\det(x)$ on the right. Then we have:

\begin{proposition}\label{formulas}
\begin{equation}\label{propf1}
\square_xD_x = p^4D_x\square_x + p^2\Delta_x + (p^{-2}+1)
\end{equation}
\begin{equation}\label{propf2}
\Delta_xD_x = p^2D_x\Delta_x + (p^{-2}+1)D_x
\end{equation}
\begin{equation}\label{propf3}
D_x X^l_{nm} = p^{2l-1} [2l] X^l_{nm}
\end{equation}
\end{proposition}

\begin{proof}
Using (\ref{del-det}), we obtain:
\begin{eqnarray*}
\partial_{11'} \partial_{22'} (f\det(x)) & = &
p^4 (\partial_{11'} \partial_{22'}f) \det(x) + p^{-1}qf + pq^{-1}( \partial_{22'}f)x_{22'} + \\
& & + pq ( \partial_{11'}f)x_{11'} + (pq-1)( \partial_{12'}f)x_{12'} + \\
& & (pq-1)( \partial_{21'}f)x_{21'} + p^{-1}q^{-1} (pq-1)^2 ( \partial_{22'}f)x_{22'} ~,    
\end{eqnarray*}
and
$$ \partial_{21'} \partial_{12'} (f\det(x)) =
p^4 (\partial_{21'} \partial_{12'}f) \det(x) - p^{-1}q^{-1}f - pq^{-1}( \partial_{12'}f)x_{12'} - $$
$$ - pq^{-1} ( \partial_{21'}f)x_{21'} - (pq^{-1}-1)( \partial_{11'}f)x_{11'} - 
(pq^{-1}-q^{-2})( \partial_{22'}f)x_{22'}~. $$
This yields (\ref{propf1}). Again using (\ref{del-det}), we have:
\begin{eqnarray*} \sum_{kl'} (\partial_{kl'}(f\det(x)))x_{kl'} & = & 
\left( p^2 (\partial_{11'}f)\det(x) + p^{-1}q^{-1}fx_{22'}\right)x_{11'} + \\
& & + \left( p^2 q^{-2}(\partial_{12'}f)\det(x) - p^{-1}q^{-1}fx_{21'}\right)x_{12'} + \\
& & + \left( p^2q^2 (\partial_{21'}f)\det(x) - p^{-1}qfx_{12'}\right)x_{21'} + \\
& & + \left( p^2(\partial_{22'}f)\det(x) + p^{-1}qfx_{11'}\right)x_{22'} ~,
\end{eqnarray*}
which yields  (\ref{propf2}).

For the third identity, we have already computed $dX^l_{nm}$ in the proof of Proposition
\ref{q-lap-prop}. The expression for $D_xX^l_{nm}$ is obtained from the one for $dX^l_{nm}$
by replacing $dx_{kl'}$ by $x_{kl'}$.
\end{proof}

The formulas in Proposition \ref{formulas} easily imply that the elements (\ref{detX}) are the
eigenfunctions of $\Delta_x$ and $\tilde{\square}_x$. Let:
$$  \tilde{\square}_x (\det(x))^k X^l_{nm} = c_k  (\det(x))^k X^l_{nm} $$
$$  \Delta_x (\det(x))^k X^l_{nm} = d_k  (\det(x))^k X^l_{nm} $$
Then the commutation relations (\ref{propf1}) and (\ref{propf2}) yield the following recurrent
formulas:
$$ c_{k+1} = p^4c_k + p^2d_k + p^{-2} + 1 ~~ , ~~ {\rm and} ~~,~~ d_{k+1} = p^2d_k + p^{-2} + 1 $$
and we know from Propositions \ref{q-lap-prop} and \ref{formulas} that $c_0=0$ and $d_0=p^{2l-1}[2l]$.
A simple induction yields $c_k=p^{2k+2l-3}[k][k+2l+2]$, as claimed in Proposition \ref{laptil}.

Finally, we will prove Proposition \ref{oast}. We note first a well-known fact in $q$-calculus: if
$$ ab = q^{-2} ba $$
\begin{equation} \label{B8}
{\rm then} ~~~ (a+b)^n = \sum_r \left\{\begin{array}{c}n\\r\end{array}\right\} a^rb^{n-r}
\end{equation}
$$ {\rm where} ~  \left\{\begin{array}{c}n\\r\end{array}\right\} = \frac{\{n\}!}{\{r\}!\{n-r\}!} ~ , $$
$$ \{n\}\! = \{1\}...\{n\} ~, ~ {\rm and} ~ \{n\} = \frac{q^{2n}-1}{q^2-1} ~. $$ 

Now the expression (\ref{q-lap-soln}) immediately implies that, up to a constant, $X^l_{nm}$
is equal to a sum:
\begin{equation}\label{B9}
\sum_r \frac{x_{11'}^r}{\{r\}!} \frac{x_{21'}^{l-m-r}}{\{l-m-r\}!}
\frac{x_{12'}^{l-n-r}}{\{l-n-r\}!} \frac{x_{22'}^{r+m+n}}{\{r+m+n\}!}
\end{equation}
Since $\det(x) =(\det(y))^{-1}$, it is enough to show that (\ref{B9}) is equal to
$(\det(x))^{2l}Y^l_{nm}$. Expanding as above we get:
\begin{equation}\label{B10}
(\det(x))^{2l} \sum_r \frac{y_{11'}^r}{\{r\}!} \frac{y_{12'}^{l-n-r}}{\{l-n-r\}!}
\frac{y_{21'}^{l-m-r}}{\{l-m-r\}!} \frac{y_{22'}^{r+m+n}}{\{r+m+n\}!}
\end{equation}
Expressing $y_{kl'}$ as $x_{kl'}/\det(x)$, commuting all factors of $(\det(x))^{-1}$
to the left and finally commuting $x_{12'}^{l-n-r}$ and $x_{21'}^{l-m-r}$, we reduce
(\ref{B10}) to a sum of monomials as in (\ref{B9}), but with certain coefficients given
by powers of $q$. An easy calculation shows that these powers do not depend on the
summation index $r$, implying that the expressions in (\ref{B9}) and ({\ref{B10}) are
indeed proportional. This concludes the proof of Proposition \ref{oast}.


 \end{document}